\documentclass[final]{autart}

\usepackage{url}
\usepackage{amsmath}
\usepackage{amssymb}
\usepackage{makeidx}
\usepackage{multicol}
\usepackage[bottom]{footmisc}
\usepackage{subfig}
\usepackage{graphicx}
\usepackage{color}
\usepackage{natbib}
\usepackage{epstopdf}
\usepackage{algorithm}
\usepackage{algcompatible}
\usepackage{calc}
\newtheorem{approxi}{Approximation}
\newtheorem{assump}{Assumption}
\newcommand{\proof}{\noindent {\bf Proof. }}
\newcommand{\app}{\begin{approxi}}
\newcommand{\eapp}{\end{approxi}}
\newcommand{\ass}{\begin{assump}}
\newcommand{\eass}{\end{assump}}
\newcommand{\teo}{\begin{thm}}
\newcommand{\eteo}{\end{thm}}
\newcommand{\corr}{\begin{cor}}
\newcommand{\ecorr}{\end{cor}}
\newcommand{\pro}{\begin{prop}}
\newcommand{\epro}{\end{prop}}
\newcommand{\lemma}{\begin{lem}}
\newcommand{\elemma}{\end{lem}}
\newcommand{\pb}{\begin{prob}}
\newcommand{\epb}{\end{prob}}
\newcommand{\df}{\begin{defn}}
\newcommand{\edf}{\end{defn}}
\newcommand{\rema}{\begin{rem}}
\newcommand{\erema}{\end{rem}}

\newcommand{\al}[1]{\begin{align} #1 \end{align}}
\newcommand{\nn}{\nonumber}

\newcommand{\Sp}[2]{\left< #1,#2 \right> }

\newcommand{\tr}{\mathop{\rm tr}}  

\newcommand{\Cc}{ \mathcal{C}}

\newcommand{\Ec}{ \mathcal{E}}

\newcommand{\Gc}{ \mathcal{G}}

\newcommand{\Ic}{ \mathcal{I}}

\newcommand{\Lc}{ \mathcal{L}}
\newcommand{\Mc}{ \mathcal{M}}

\newcommand{\Qc}{ \mathcal{Q}}

\newcommand{\Sc}{ \mathcal{S}}

\newcommand{\Vc}{ \mathcal{V}}

\newcommand{\Cs}{ \mathbb{C}}

\newcommand{\Es}{ \mathbb{E}}

\newcommand{\Ns}{ \mathbb{N}}

\newcommand{\Rs}{ \mathbb{R}}

\newcommand{\Zs}{ \mathbb{Z}}

\newcommand{\Dd}  {\mathrm{D}}
\newcommand{\Ld}  {\mathrm{L}}

\newcommand{\Td}  {\mathrm{T}}

\newcommand{\Mb}  {\mathbf{M}}
\newcommand{\Qb}  {\mathbf{Q}}

\begin{document}

\begin{frontmatter}
\title{Empirical Bayesian Learning in AR Graphical Models}


\author[Padova]{Mattia Zorzi}
\ead{zorzimat@dei.unipd.it}

\address[Padova]{Dipartimento di Ingegneria dell'Informazione, Universit\`a degli studi di
Padova, via Gradenigo 6/B, 35131 Padova, Italy}


\begin{abstract} 
We address the problem of learning graphical models which correspond to  high dimensional autoregressive stationary stochastic processes. A graphical model describes the conditional dependence relations among the components of a stochastic process and represents an important tool in many fields. We propose an empirical Bayes estimator of sparse autoregressive graphical models and latent-variable autoregressive graphical models. Numerical experiments show the benefit to take this Bayesian perspective for learning these types of graphical models.
\end{abstract}
\end{frontmatter}

\section{Introduction}
In modern applications many variables are accessible to observation. In some cases the latter can be modeled with a high dimensional Gaussian random vector.
To gain some insight about the relation among those variables we can attach to it a graphical model \citep{LAURITZEN_1996,Willsky02multiresolutionmarkov}. The latter is an undirected graph wherein nodes correspond to the components (i.e. variables) of the random vector and there is the lack of an edge between two nodes if and only if the corresponding variables are conditionally independent given the others. It turns out that sparse graphical models, i.e. graphs with few edges, have the inverse covariance matrix of the random vector which is sparse. Then, the problem of estimating a sparse graphical model from the observed data can be formulated as a regularization problem: find such an inverse matrix which minimizes the negative log-likelihood and a regularization term inducing sparsity \citep{banerjee2008model}.

An important aspect is that variables are typically measured over time and can thus be modeled as a high dimensional autoregressive (AR) Gaussian stationary stochastic process. Then, we can attach a graphical model describing the conditional dependence relations among the variables. It is possible to prove that sparse graphical models have the inverse power spectral density (PSD) of the process which is sparse. \cite{SONGSIRI_TOP_SEL_2010} proposed a regularized estimator for estimating sparse AR graphical models in the same spirit of \citep{banerjee2008model}. 
\cite{ARMA_GRAPH_AVVENTI} showed that the aforementioned estimator is a relaxed version of a maximum entropy estimator. The latter solves a covariance extension problem whose dual problem does coincide with the one proposed in \cite{SONGSIRI_GRAPH_MODEL_2010}. As a consequence, 
the estimator proposed by \cite{SONGSIRI_TOP_SEL_2010} is strictly connected with the generalized moment problems in the sense
of Byrnes-Georgiou-Lindquist which have been extensively studied by many researchers, e.g. \cite{A_NEW_APPROACH_BYRNES_2000,FERRANTE_TIME_AND_SPECTRAL_2012,karlsson2013uncertainty,DUAL,BETA,6909016,ALPHA}.  Since then, many other extensions has been proposed: \cite{MAANAN2017122} proposed a two stage approach to estimate sparse AR graphical models; \cite{REC_SPARSE_GM} proposed a regularized estimator for sparse graphical models of reciprocal processes; \cite{Chandrasekaran_latentvariable}, \cite{LATENTG}, \cite{e20010076}, \cite{CDC_BRAIN15} and \cite{CICCONE}
  proposed regularized estimators for the so called latent-variable graphical models.

The regularizers for inducing sparsity in such graphical models are $\ell_1$-like norms. The $\ell_1$ norm, however, penalizes differently the nonnull coefficients: larger coefficients are penalized more heavily than smaller coefficients. This imbalance produces an estimator with a remarkable mean squared error.
\cite{candes2008enhancing} proposed to reduce this imbalance by considering a weighted $\ell_1$ norm wherein the weights are computed in an iterative fashion. The resulting procedure is similar to the adaptive lasso and is typically called iterative reweighted algorithm \citep{WIPF_2010}.

The present paper proposes a regularized estimator in the spirit of \cite{SONGSIRI_TOP_SEL_2010} for sparse AR graphical models where the $\ell_1$-like norm is substituted by a weighted $\ell_1$-like norm leading to an iterative reweighted procedure. Interestingly, drawing inspiration by \cite{asadi2009map, scheinberg2010sparse},
the proposed method can be understood as an empirical Bayes approach which provides a suitable updating rule for the weights. Such idea is then extended to the regularized estimator for latent-variable AR graphical models proposed in \cite{LATENTG}.

The outline of the paper is as follows. In Section \ref{sec:AR_GM} we introduce
the problem of estimating sparse AR graphical models. 
In Section \ref{sec_reweight} we propose an iterative reweighted for solving such a problem, while in Section \ref{sec:Bayes_S} we derive the estimator using a Bayesian perspective. Section \ref{sec:Bayes_SL} regards the identification of latent-variable AR graphical models. Section \ref{sec:sim} contains some numerical experiments to test the performance of the proposed estimators. Finally, the conclusions are drawn in Section \ref{sec:concl}.

\subsection*{Notation}
The vector space $\Rs^m$ is endowed with the inner product $\Sp{x}{y}=x^T y$. $x\geq 0$ ($x>0$) 
with $x\in \Rs^m$ means that all the entries of the vector are nonnegative (positive). 
The vector space $\Rs^{m\times m}$ is endowed with the inner product $\Sp{X}{L}=\tr(X L^T)$. $\Qb_{m}$ denotes the vector space of symmetric matrices of dimension $m\times m$, if $X\in\Qb_m$ is
positive definite (semi-definite) we write $X\succ 0$ ($X\succeq 0$). $|X|$ denotes the determinant of matrix $X\in \Qb_m$. $\mathrm{diag}(X)$ is the diagonal matrix whose main diagonal coincides with the one of $X$. $(X)_{jh}$ denotes the entry in position $(j,h)$ of matrix $X$. A matrix $A\in\Rs^{l\times m(n+1)}$ with $l\leq m$ will be partitioned as
$A=[\, A_0 \, A_1 \, \ldots \, A_n \,] $ with $A_j\in\Rs^{l\times m}$. $\Mb_{m,n}$ is the vector space of matrices $Y:=[\, Y_0 \, Y_1 \, \ldots \, Y_{n}\,]$ with  $Y_0\in\Qb_m$ and $Y_1\ldots Y_{n}\in\Rs^{m\times m}$. The corresponding inner product is $\Sp{Y}{Z}=\tr(YZ^T)$. The linear operator $\Td: \Mb_{m,n} \rightarrow \Qb_{m(n+1)}$ constructs a symmetric {\em Toeplitz}
matrix $ \Td(Y)$ where its first block row is $[\,
                      Y_{0} \, Y_{1} \, \ldots \,Y_{n}\,].$
 The adjoint operator of $\Td$ is denoted by $\Dd: \Qb_{m(n+1)}\rightarrow \Mb_{m,n}$ and defined as follows. If $X \in\Qb_{m(n+1)}$ is partitioned as a $n+1\times n+1$ block matrix 
with $X_{hj}$, $j,h=0\ldots n$, the block in position $(h,j)$
then $\Dd(X)=[\,
                        \Dd_0(X) \, \ldots \, \Dd_{n}(X) \,]
$ where $\Dd_0(X)=\sum_{h=0}^{n} X_{h h}$, $ \Dd_k(X)=2\sum_{h=0}^{n-k} X_{h\; h+k}$,  $k=1\ldots n$.
We define the index set $\Ec \subseteq \Vc\times \Vc$ with $\Vc:=\{1,2,\ldots m\}$. Functions on the unit circle $\{e^{i\vartheta} \hbox{ s.t. } \vartheta \in[-\pi,\pi]\}$ will be denoted by capital Greek letters, e.g. $\Phi(e^{i\vartheta})$ with $\vartheta\in[-\pi,\pi]$, and the dependence upon $\vartheta$ will be dropped if not needed, e.g. $\Phi$ instead of $\Phi(e^{i\vartheta})$. $\Ld_2^{m\times m}$ denotes the space of $\Cs^{m\times m}$-valued functions defined on the unit circle which are square integrable. Given $\Phi\in\Ld^{m\times m}_2$, the shorthand notation $\int \Phi$ denotes the integration of $\Phi$ taking place on the unit circle with respect to the normalized {\em Lebesgue} measure. Then, the inner product in $\Ld_2^{m\times m}$ is $\Sp{\Phi}{\Sigma}=\tr\int \Phi\Sigma^*$. Given an analytic function $\Lambda\in\Ld_2^{m\times m}$, its (normal) rank is denoted by $\mathrm{rank}(\Lambda)$.
If $\Phi(e^{i\vartheta})$ is positive definite (semi-definite) for each $\vartheta\in[-\pi,\pi]$, we will write $\Phi\succ 0$ ($\Phi \succeq 0$).
We define the following family of matrix pseudo-polynomials
  \al{ \label{set_Qc_mn_reparametrized}\Qc_{m,n}=\{\Delta X \Delta^* \hbox{ s.t. } X\in\Qb_{m(n+1)}\}}
where 
$ \Delta(e^{i\vartheta}):=[\,
                                                                            I_m \, e^{i\vartheta}I_m \, \ldots \, e^{in \vartheta} I_m \,]$
 is the {\em shift operator}. Given a $m$-dimensional stochastic process $y=\{\,y(t) ,\; t\in \Zs\,\}$, $y_j$ denotes the $j$-th component (i.e. variable) of $y$. With some abuse of notation, $y(t)$ will both denote a random vector and its sample value.
Given a function $f(u,v)$, $\nabla_u f(\bar u,\bar v)$ and $\nabla_v f(\bar u,\bar v)$ denote the gradient of $f$ with respect to $u$ and $v$, respectively, computed at $(\bar u,\bar v)$.

\section{Identification of Sparse AR Graphical Models}\label{sec:AR_GM}

\begin{figure}
\centering
\subfloat{\includegraphics[width=0.44\columnwidth]{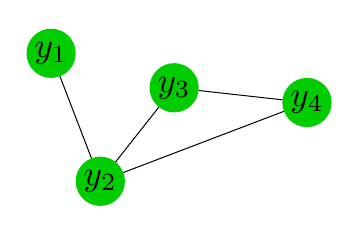} }
\subfloat{\includegraphics[width=0.54\columnwidth]{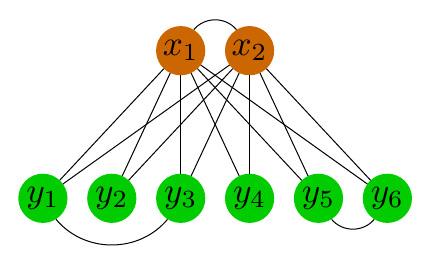}}
\caption{{\em Left.} Example of a graphical model for $y = [\, y_1\, y_2 \, y_3\, y_4 \,]^T$. {\em Right.} Example of a latent-variable graphical model: $y_1\ldots y_6$ are the manifest variables $x_1,x_2$ are the latent variables.} \label{fig:graph_ex}
\end{figure}

 Assume to collect the data $ y^N:=\{y(1), y(2) \ldots y (N)\}$ generated by the AR Gaussian discrete-time zero mean full rank stationary stochastic process $y=\{\, y(t),\;  t\in \mathbb{Z}\}$ defined as 
\al{ \label{mod_AR}y(t)=-\sum_{k=1}^{ n }A_k y(t-k) +e(t).} $y(t)$ takes values in $\Rs^m$, $A_k\in\Rs^{m\times m}$ and $e(t)$
is white Gaussian noise with covariance matrix $R\succ0$.  Both $A_k$ and $R$ are unknown. The order of the AR process, i.e. $n$, is assumed to be known. We want to estimate $A_k$ and $R$ using $y^N$. It is well known that an equivalent description of $y$ is given by its PSD
\al{\Phi(e^{i\vartheta})=\sum_{k\in \mathbb{Z}}  e^{-ik \vartheta} R_k, \; \; \; \vartheta\in[-\pi,\pi]}   where $R_k=\Es[y(t+k) y(t)^T]$, with $k\in \mathbb{Z}$, is the covariance lags sequence. Accordingly, the aforementioned problem is equivalent to estimate $\Phi$ from $y^N$.  There are situations (e.g. $m$ is large) in which one is interested to estimate a PSD reflecting only the most important conditional dependence relations among the variables of $y$. Let $I\subset \Vc$ be an arbitrary index set. We denote as 
 $\chi_I=\overline{\mathrm{span}}\{\, y_j(t) \hbox{ s.t. }\, j\in I\, , \, t\in \Zs\, \} $
the closure of the vector space of all finite linear combinations (with real coefficients)
of $y_j(t)$ with $j\in I $ and $t\in \Zs$. Let $j\neq h$, we say that $y_j$ and $y_h$ are conditionally independent 
given the other variables if $\chi_{\{ j\}} \, \bot\, \chi_{\{h\}} \, |\, \chi_{\Vc \setminus \{j,h\}}.$
These conditional dependence relations define an interaction graph $\Gc(\Vc,\Ec)$ where $\Vc$ and $\Ec$ denote the set of nodes and edges, respectively. More precisely, the nodes represent the variables $y_1,y_2\ldots y_m$ and the lack of an edge means conditional independence \citep{REMARKS_BRILLINGER_1996}:
\al{(j,h) \notin \Ec \, \, \iff \, \, \chi_{\{ j\}} \, \bot\, \chi_{\{h\}} \, |\, \chi_{\Vc \setminus \{j,h\}}, \; \;  j\neq h.}
An example of graphical model is provided in Figure \ref{fig:graph_ex} (left).
 \cite{Dahlhaus2000} proved that $y_j$ and $y_h$ are conditionally independent if and only if $(\Phi^{-1}(e^{i\vartheta}))_{jh}=0$, $\forall \, \vartheta\in [-\pi,\pi].$ This characterization allows to infer conditional independence relations by promoting sparsity 
in the estimation of the inverse PSD for model (\ref{mod_AR}). Since $y$ is an AR process of order $n$, we can parametrize its PSD as $\Phi=\Sigma^{-1}$ where $\Sigma=S_0+\frac{1}{2}\sum_{k=1}^n  S_k e^{-ik\vartheta}+S_k^T e^{ik\vartheta}\in \Qc_{m,n}$
and we define $S=[\, S_0\, S_1\,\ldots \, S_n \,]\in \Mb_{m,n}$. Then, a regularized maximum likelihood (ML) estimator of $\Sigma$, and thus of $\Phi$, is given by solving \citep{SONGSIRI_TOP_SEL_2010}:
\al{\label{pb_td1}\hat \Sigma= &\underset{\Sigma\in \Qc_{m,n}}{ \mathrm{argmin}} \,\ell( y^N; \Sigma)+ \gamma\bar h_{\infty}(\Sigma)\nn\\
& \hbox{ s.t. }  \Sigma\succ0 .}
The term $\ell( y^N; \Sigma)$ is an approximation of the negative log-likelihood of $  y(n
+1) \ldots  y(N)$ given $  y(1) \ldots  y (n)$ under model (\ref{mod_AR}): 
\al{\label{cond_PDF}\ell(  y^N; \Sigma)= -\frac{N-n}{2}\int \log |\Sigma| +\frac{N-n}{2}\Sp{\Phi_{ y^N}}{\Sigma}+c} 
where \al{\label{windowed_cor}& \Phi_{y^N}=\hat R_0+\frac{1}{2}\sum_{k=1}^n  \hat R_k e^{-ik\vartheta}+\hat R_k^T e^{ik\vartheta}\nn \\
& \hat R=[\, \hat R_0\,\ldots  \, \hat R_n \,], \; 
\hat R_k =\frac{1}{N-n} \sum_{t=1}^{N-k} y(t+k)y(t)^T
} and $c$ is a term not depending on $y^N$ and $\Sigma$. It is worth noting that $\hat R_k$ represents an estimate of $R_k$ from data $ y^N$, indeed  $\Phi_{y^N}$ is the truncated periodogram of $\Phi$ computed from $y^N$. Under the assumption that $y$ is a full rank process, then for $N$ sufficiently large we have that $\Td(\hat R)\succ 0$ with high probability. Accordingly, throughout the paper we make the assumtion that $\Td(\hat R)\succ 0$. The penalty term  $\bar h_\infty\, :\, \Qc_{m,n} \longrightarrow \Rs$ \al{\label{h_inf} \bar h_\infty &(\Sigma) =\sum_{j>h}  q_{jh}(\Sigma)}  
with \al{ q_{jh}(\Sigma):= \max\{ |(S_0)_{jh}|, \max_{k=1 \ldots n} |(S_k)_{jh}|, \max_{k=1 \ldots n} |(S_k)_{hj}|\}\nn}
encourages a common sparsity pattern (i.e. group sparsity) on the coefficients $S_k$ of $\Sigma$. Finally, $\gamma>0$ is the regularization parameter. Since $\Sigma\in \Qc_{m,n}$ and in view of (\ref{set_Qc_mn_reparametrized}), Problem (\ref{pb_td1}) can be rewritten in terms of a matrix $X\in\Qb_{m(n+1)}$ \citep{SONGSIRI_TOP_SEL_2010}:\al{\label{pb_td2}\hat X= &\underset{X\in \Qb_{m(n+1)}}{ \mathrm{argmin}} \ell (y^N; X)+\gamma h_{\infty}(\Dd(X))\nn\\ & \hbox{ s.t. }  X_{00}\succ 0,\; \;X\succeq 0}
where $\ell ( y^N; X)=-\frac{N-n}{2}\log |X_{00}|+\frac{N-n}{2}\tr(\Td(\hat R) X)+c$; function  $h_\infty\, :\, \Mb_{m,n} \longrightarrow \Rs$ is defined as (\ref{h_inf}) but now the domain is replaced by $\Mb_{m,n}$. Finally, the optimal solution $\hat \Phi= \hat \Sigma^{-1}$ is such that $\hat \Sigma=\Delta\hat X\Delta^*$. It is clear that  $\hat \Phi^{-1}$ depends on $\gamma$. \cite{SONGSIRI_TOP_SEL_2010} 
proposed to select $\gamma$ by computing $j$ values of $\gamma$ according to the so called ``trade-off'' curve. Then the corresponding candidate models are ranked by using a BIC criterium. It is worth noting that the latter applies a thresholding on the partial coherence of the estimated PSD in order to measure the complexity (in terms of number of edges) of the candidate model.
 
\section{Iterative Reweighted Method} \label{sec_reweight}
In this section we investigate the possibility to modify the penalty term $\bar h_\infty$ in (\ref{pb_td1}) 
in such a way to improve the ability to estimate the support of $\Phi^{-1}$.  Notice that $h_\infty$ can be understood as the $\ell_1$ norm of the vector $[\, q_{11}(\Sigma) \,q_{12}(\Sigma) \ldots q_{nn}(\Sigma) \,]^T$ which represents the convex surrogate of the corresponding $\ell_0$ norm. As highlighted in \cite{candes2008enhancing}: ``A key difference between the $\ell_0$ and the $\ell_1$ norms is the dependence on magnitude: larger coefficients are penalized more heavily in the $\ell_1$ norm than smaller coefficients, unlike the more democratic penalization of the $\ell_0$ norm''. We address this imbalance by considering a weighted penalty function $\bar h_{W}(\Sigma) =\sum_{j \geq h}  \gamma_{jh} q_{jh}(\Sigma)$ where $\Sigma\in \Qc_{m,n}$ 
and $\gamma_{jh}> 0$. Note that, we penalize also the entries in the main diagonal  of $\Sigma$. The latter regularization is not imposed in order to obtain sparsity in the main diagonal of $\Sigma$, otherwise constraint $\Sigma\succ 0$ is no longer satisfied, but rather it is imposed in order to reduce the variance of the estimator for the variables in the main diagonal of $\Sigma$. Accordingly, we will expect to find $\gamma_{jj}>0$ but smaller than $\gamma_{jh}$ with $j\neq h$. Thus, we modify Problem (\ref{pb_td1}) as
\al{\label{pb_rw1}\hat \Sigma= &\underset{\Sigma \in \Qc_{m,n}}{ \mathrm{argmin}} \, \ell( y^N; \Sigma)+ \bar h_{W}(\Sigma)\nn\\
&\hbox{ s.t. }  \Sigma\succ0.}
Following the same reasoning in \cite{SONGSIRI_TOP_SEL_2010} we can reformulate Problem (\ref{pb_rw1}) as follows:
\al{\label{pb_rw2}\hat X= &\underset{X\in \Qb_{m(n+1)}}{ \mathrm{argmin}} \ell ( y^N; X)+\bar h_{W}(\Dd(X))\nn\\ & \hbox{ s.t. }  X_{00}\succ 0, \; \;X\succeq 0} where $h_W\, : \, \Mb_{m,n} \rightarrow \Rs$ is defined as $\bar h_W$ but now the domain is $\Mb_{m,n}$. It is worth noting that (\ref{pb_rw1}) and (\ref{pb_rw2}) are equivalent provided that the optimal solution $\hat X$ is such that $\Delta \hat X \Delta^*\succ 0$.

\teo \label{teo_S1}The dual problem of (\ref{pb_rw2}) is  
\al{\label{dual_rw2} &\underset{\substack{Z\in \Mb_{m,n}\\ W\in \Qb_{m}}}{ \mathrm{max}} \log | W|+m\nn\\ & \hbox{ s.t. } W\succ 0\nn\\& \hspace*{0.7cm}\Td(\hat R)+\Td(Z)\succeq \left[\begin{array}{cc}W &  0 \\0  & 
0\end{array}\right]\nn\\ & \hspace*{0.7cm} \sum_{k=0}^n |(Z_k)_{jh}|+|(Z_k)_{hj}|\leq \frac{2\gamma_{jh}}{N-n},\; \; j>h\nn\\
 &\hspace*{0.7cm} \sum_{k=0}^n |(Z_k)_{jj}|\leq \frac{2\gamma_{jj}}{N-n},\; \; j=h.} The latter does admit a unique solution $\hat X$ such that $\Delta \hat X \Delta^*\succ 0$. Accordingly (\ref{pb_rw1}) and (\ref{pb_rw2}) are equivalent.\eteo
Theorem \ref{teo_S1} is important not only for establishing the existence of the solution of (\ref{pb_rw1}) and (\ref{pb_rw2}), but it provides also a tractable formulation for the computation of the solution, see Section \ref{sec:imple}.

It remains to select a suitable set of weights $\gamma_{jh}$. The latter should counteract the magnitude imbalance which characterizes the $\ell_1$-like norm; more precisely, $\gamma_{jh}$ should be inversely proportional to $q_{jh}(\Sigma)$. For instance, we could take $\gamma_{jh}=(q_{jh}(\Sigma)+\varepsilon_S)^{-1}$ with $\varepsilon_S>0$ sufficiently small. However, $\Sigma$ is unknown. Accordingly, we propose an iterative reweighted algorithm \citep{WIPF_2010} which constructs iteratively a set of weights by using the information from the current estimate of $\Sigma$ (i.e. $X$), see Algorithm \ref{algo:RWS}
\begin{algorithm}
\caption{Iterative Reweighted Algorithm}
\label{algo:RWS}
\begin{algorithmic}[1] \small
\STATE $l=0$
\STATE Initialize $\hat \gamma_{jh}^{(0)}$ with $j,h=1\ldots m$
\REPEAT 
\STATE Solve Problem (\ref{pb_rw2}) with $\gamma_{jh}=\hat{\gamma}_{jh}^{(l)}$; let $\hat X^{(l)}$ denote 
\STATEx \hspace{0.3cm}the corresponding solution
\STATE $q_{jh}(\Dd(\hat X^{(l)}))\leftarrow\max\{ |(\Dd_0(\hat X^{(l)}))_{jh}|,$ \STATEx  \hspace{1.2cm} $\max_{k>0} |(\Dd_k(\hat X^{(l)}))_{jh}|, \max_{k>0} |(\Dd_k(\hat X^{(l)}))_{hj}|\}$
\STATE Update the weights $\hat \gamma_{jh}^{(l+1)}\leftarrow (q_{jh}(\Dd(\hat X^{(l)}))+\varepsilon_S)^{-1}$
\STATE $l\leftarrow l+1$
\UNTIL{ $\|\hat X^{(l)}-\hat X^{(l-1)}\|\geq \varepsilon $ \textbf{or} $l\leq l_{\mathrm{MAX}}$
\STATE $\hat \Phi=(\Delta \hat X^{(l)}\Delta^*)^{-1}$ }
\end{algorithmic}
\end{algorithm}
where we recall that $S_k=\Dd_k(X)$. Parameter $\varepsilon_S>0$ in Step 6 ensures that a zero-valued entry in positions $(j,h)$ and $(h,j)$ does not strictly prohibit a nonnull estimate at the next step. Although Algorithm \ref{algo:RWS} has been introduced by an heuristic reasoning, it can be interpreted as the Majorization-Minimization (MM) algorithm for the problem  \al{\label{pb_MM} &\underset{\substack{X\in \Qb_{m(n+1)}}}{ \mathrm{argmin}}\; \; \sum_{j\geq h}\log ( q_{jh}(\Dd(X))+\varepsilon_S )+\ell ( y^N; X)\nn\\ & \hbox{ s.t. }  X_{00}\succ 0,\; \; X\succeq 0.}   The latter provides a regularized ML estimator of $\Sigma$ through the relation $S_k=\Dd_k(X)$ and $S_k$ are the coefficients of $\Sigma \in\Qc_{m,n}$. The log-sum penalty induces sparsity on $\Sigma$. It is well known that such a penalty outperforms the corresponding $\ell_1$-like norm for estimating the correct sparsity pattern, see \cite{candes2008enhancing}. On the other hand, the log-sum penalty is concave making Problem (\ref{pb_MM}) nonconvex. In order to see the aforementioned connection, we can rewrite 
(\ref{pb_MM}) as 
\al{\label{pb_MM2}&\underset{\substack{X\in \Qb_{m(n+1)}\\ u_{jh}\in\Rs,\; j\geq h}}{ \mathrm{argmin}} \; \;\sum_{j\geq h}\log (u_{jh}+\varepsilon_S )+\ell ( y^N; X)\nn\\ & \hbox{ s.t. }  X_{00}\succ 0,\; \; X\succeq 0\nn\\
& \hspace{0.7cm}u_{jh}\geq q_{jh}(\Dd(X)).} Indeed, if $\hat X$ is the optimal solution of (\ref{pb_MM}), then the optimal solution of (\ref{pb_MM2}) is $\hat X$ and  $\hat u_{jh}=q_{jh}(\Dd(\hat X))$ with $j\geq h$. Let $u\in\Rs^{m(m+1)/2}$ be the vector obtained by stacking $u_{jh}$ with $j\geq h$. We define \al{f(u,X)=\sum_{j\geq h}\log (u_{jh}+\varepsilon_S )+\ell ( y^N; X).\nn} 
The latter is majorized by:
\al{g(u&,X;\hat u^{(l)},\hat  X^{(l)})\nn\\ &:=f(\hat u^{(l)},\hat X^{(l)})+\Sp{\nabla_u f(\hat u^{(l)})}{u-\hat u^{(l)}}+\ell ( y^N; X)\nn\\ 
&=f(\hat u^{(l)}, \hat X^{(l)})+\sum_{j\geq h} \frac{u_{jh}-\hat u_{jh}^{(l)}}{\hat u_{jh}^{(l)}+\varepsilon_S}+\ell ( y^N; X).\nn} 
More precisely, the first term in $f$, which is a concave function in $u$, is majorized by the tangent at $\hat u^{(l)}$.
Then, the MM algorithm is
\al{\label{pb_MM3uno}(\hat u^{(l+1)},\hat X^{(l+1)})= & \underset{u,X}{ \mathrm{argmin}} \; \; g(u,X;\hat u^{(l)},\hat X^{(l)})\nn\\ & \hbox{ s.t. }  X_{00}\succ 0,\; \; X\succeq 0\nn\\
& \hspace{0.7cm}u_{jh}\geq q_{jh}(\Dd(X)).} Removing the terms in $g$ not depending on $u$ and $X$, we obtain
\al{\label{pb_MM3}(u^{(l+1)},X^{(l+1)})= & \underset{u,X}{ \mathrm{argmin}} \; \; \sum_{j\geq h} \frac{u_{jh}}{\hat u_{jh}^{(l)}+\varepsilon_S}+\ell ( y^N; X)
\nn\\ & \hbox{ s.t. }  X_{00}\succ 0,\; \; X\succeq 0\nn\\
& \hspace{0.7cm}u_{jh}\geq q_{jh}(\Dd(X)).}
Since we have $u_{jh}=q_{jh}(X)$ for the optimal solution of (\ref{pb_MM3}), then the latter problem is equivalent to
\al{\label{pb_MM4}\hat X^{(l+1)}= & \underset{X}{ \mathrm{argmin}} \; \; \sum_{j\geq h} \frac{q_{jh}(X)}{q_{jh}(\hat X^{(l)})+\varepsilon_S}+\ell ( y^N; X)
\nn\\ & \hbox{ s.t. }  X_{00}\succ 0,\; \; X\succeq 0.} 
 Defining \al{ \label{com_gamma}\hat \gamma_{jh}^{(l+1)}:=\frac{1}{q_{jh}(\Dd(\hat X^{(l)}))+\varepsilon_S} } and substituting it in (\ref{pb_MM4}) we obtain (\ref{pb_rw2}) where $\gamma_{jh}$ has been replaced by $\hat\gamma_{jh}^{(l+1)}$. In other words, (\ref{pb_MM4}) is equivalent to: compute $\hat \gamma_{jh}^{(l+1)}$ as in (\ref{com_gamma}) and then solve (\ref{pb_rw2}) with $\gamma_{jh}=\hat \gamma_{jh}^{(l+1)}$. The latter procedure coincides with the iterative reweighted scheme in Algorithm \ref{algo:RWS}. 

\pro \label{convMM}Consider the sequence $\hat X^{(l)}$ generated by Algorithm \ref{algo:RWS}. Then, $\hat X^{(l)}$ converges to the set of stationary points of Problem (\ref{pb_MM}).   
\epro

Simulation evidence showed that the proposed iterative reweighted algorithm does not perform well, that is the estimated PSD is not close to the actual one, and neither the sparsity pattern of the inverse. 
Even if the updating rule for $\gamma_{jh}$s seems reasonable, it is not the best choice that
we can apply. In the next section we will see how to find a better update.

\section{A Bayesian Perspective}\label{sec:Bayes_S}
Till now the problem of estimating $\Sigma$ sparse has been considered  according to the Fisherian perspective that is $\Sigma$ is an unknown but fixed function in $\Qc_{m,n}$, i.e. the parameters $S_k$, $k=0\ldots n$, characterizing $\Sigma$ are unknown but fixed quantities. In this section we   propose a method which is based on a Bayesian perspective that is $\Sigma$ is a stochastic process taking values in $\Qc_{m,n}^+:=\{\Sigma \hbox{ s.t. } \Sigma\succ 0\}$. This means that the parameters characterizing $\Sigma$ are random variables with a suitable PDF or simply prior. Let $p(\Sigma)$ be the prior of $\Sigma$. We recall that $(\Sigma)_{jh}$ is the entry of $\Sigma$ in position $(j,h)$:
 \al{(\Sigma(e^{i\vartheta}))_{jh}=(S_0)_{jh}+\frac{1}{2}\sum_{k=1}^n (S_k)_{jh}e^{-i k \vartheta}+(S_k)_{hj}e^{i k \vartheta}.\nn}
 We assume that $(\Sigma)_{jh}$s are independent each other, accordingly 
\al{p_{\pmb{\gamma}}(\Sigma)=\prod_{j\geq h} p_{\gamma_{jh}}((\Sigma)_{jh}),\nn}
and 
\al{ p_{\gamma_{jh}}&((\Sigma)_{jh})=\frac{e^{-\gamma_{jh}q_{jh}(\Sigma)}}{c_{jh}}\nn }
where $c_{jh}$ is the normalizing constant and $\gamma_{jh}> 0$, $j\geq h$, are referred to as hyperparameters. 
The latter are the parameters characterizing the prior. Therefore, we have
\al{\label{p_Sigma}p_{\pmb \gamma}(\Sigma)=\frac{e^{-\bar h_W(\Sigma)}}{\prod_{j\geq h} c_{jh}}}
where $\pmb{\gamma}\in\Rs^{m(m+1)/2}$ denotes the hyperparameters vector containing $\gamma_{jh}$ with $j\geq h$.
The negative log-likelihood of $y^N$ and $\Sigma$ takes the form:
\al{\ell ( y^N,\Sigma;\pmb \gamma)&=-\log p_{\pmb \gamma}( y^N,\Sigma)\nn\\
&=-\log p( y^N|\Sigma)-\log p_{\pmb \gamma}(\Sigma).}
It is clear that the negative log-conditional PDF $-\log p( y^N|\Sigma)$ does coincide with (\ref{cond_PDF}), thus 
\al{\ell ( y^N,\Sigma;\pmb \gamma)=\ell( y^N; \Sigma)+\bar h_W(\Sigma)+\sum_{j\geq h}\log c_{jh} }
where the last term does not depend on $\Sigma$.
We conclude that the MAP estimator  of $\Sigma$ is given by (\ref{pb_rw1}). This result is not surprising, indeed it is well known that MAP estimators in a Bayesian perspective can be interpreted as regularized estimators in the Fisherian perspective. The substantial difference between the two perspectives is that the former provides the way to estimate the hyperparameters vector $\pmb \gamma$. We define as the negative log-marginal likelihood: 
 \al{\label{ML_S}\ell(y^N; \pmb\gamma):= -\log \int _{\Qc_{m,n}^+} p_{\pmb \gamma}( y^N, \Sigma ) \mathrm d \Sigma.} An estimate of $\pmb \gamma$ is given by the empirical Bayes approach \citep{friedman2001elements}:
\al{\label{emp_Bayes}\hat{\pmb{ \gamma}}= &\underset{{\pmb \gamma}> 0}{ \mathrm{argmin}}\; \ell( y^N; \pmb \gamma).}  Then, the MAP estimator of $\Sigma$ is given by (\ref{pb_rw1}) with $\pmb \gamma =\hat{\pmb \gamma}$. However, it is not possible to find an analytical expression for $ \ell( y^N; \pmb \gamma)$ making challenging the optimization of $\pmb \gamma$.

 \pro \label{prop_norm_const}Consider the prior of $\Sigma$ defined in (\ref{p_Sigma}). Then, we have 
\al{c_{jh}\leq \left\{\begin{array}{ll} \upsilon_{jh}\gamma_{jh}^{-(n+1)}, &  \hbox{ if $j=h$}; \\ \upsilon_{jh} \gamma_{jh}^{-(2n+1)}, & \hbox{ if $j>h$}\end{array}\right.}
where the terms $\upsilon_{jh}$ in the relation above do not depend on $\pmb \gamma$.
\epro
Accordingly,
\al{\label{lik_proportional}\ell( &y^N , \Sigma; \pmb \gamma)\leq \tilde \ell( y^N , \Sigma; \pmb \gamma):= \ell( y^N; \Sigma)+\bar h_W(\Sigma)\nn\\ &-\sum_{j>h} (2n+1)\log \gamma_{jh} -\sum_{j=1}^n (n+1) \log \gamma_{jj}+c} 
where $c$ is a term not depending on $ y^N$, $\Sigma$ and $\pmb \gamma$. An alternative simplified approach to estimate $ \pmb \gamma$ is the generalized maximum likelihood (GML) method, \cite{zhou1997approximate}: instead of maximizing $\ell ( y^N; \pmb \gamma)$ with respect to $\pmb \gamma$, the latter
is computed with $\Sigma$ as the pair $(\hat{\pmb \gamma},\hat \Sigma)$ that jointly minimizes $  \tilde \ell ( y^N,\Sigma;\pmb \gamma)$. Then, the optimization can be performed in a two-step algorithm:
\al{\label{2step_sigma}\hat \Sigma^{(l)} =& \underset{  \Sigma}{\;\mathrm{argmin}}\;  \tilde \ell ( y^N,\Sigma;\hat{\pmb \gamma}^{(l)})\\ 
\label{2step_gamma}\hat{\pmb \gamma}^{(l+1)} =& \underset{\pmb\gamma }{\;\mathrm{argmin}}\; \tilde \ell ( y^N,\hat \Sigma^{(l)};\pmb \gamma).}
Step (\ref{2step_sigma}) is the MAP estimator of $\Sigma$ given the current choice of $\pmb \gamma$ which is equivalent to (\ref{pb_rw1}). Step (\ref{2step_gamma}) is the estimator of $\pmb \gamma$ using the current MAP estimate of $\Sigma$ as a direct observation. Notice that (\ref{2step_gamma}) is equivalent to 
\al{\hat{\gamma}^{(l+1)}_{jh} =  \left\{\begin{array}{ll} \; \underset{\gamma_{jh}> 0 }{\mathrm{ argmin}}  \;\gamma_{jh} q_{jh}(\hat S^{(l)})-(n+1)\log \gamma_{jh}, &     \\ 
& \hspace{-1.3cm}\hbox{if } j=h;\\ \; \underset{\gamma_{jh}>0 }{\mathrm{argmin}}\; \gamma_{jh} q_{jh}(\hat S^{(l)})-(2n+1)\log \gamma_{jh}, &  \\
&  \hspace{-1.3cm}\hbox{if } j>h  \\ \end{array}\right. \nn}
where we recall that $\hat S^{(l)}=\Dd(\hat X^{(l)})$.
\pro \label{prop_opt_gamma}Under the assumption that $q_{jh}(\hat S^{(l)})>0$, it holds that 
\al{\label{rew_MLsenza_e}\hat{\gamma}^{(l+1)}_{jh}=\left\{\begin{array}{ll} \frac{n+1}{q_{jh}(\hat S^{(l)})},& \hbox{ if } j=h  \\ \frac{2n+1}{q_{jh}(\hat S^{(l)})}, & \hbox{ if } j>h. \end{array}\right. }
\epro
To deal also with the case that $q_{jh}(\hat S^{(l)})=0$, we consider the modified updating 
\al{\label{rew_ML}\hat{\gamma}^{(l+1)}_{jh}=\left\{\begin{array}{ll} \frac{n+1}{q_{jh}( 	\hat S^{(l)})+\varepsilon_S },& \hbox{ if } j=h  \\ \frac{2n+1}{q_{jh}(\hat S^{(l)})+\varepsilon_S }, & \hbox{ if } j>h \end{array}\right.  } with $\varepsilon_S>0$. It is not difficult to see that the latter modification is equivalent to assume that $\gamma_{jh}$s are modeled as independent random variables with exponential hyperprior: $p(\gamma_{jh})=\varepsilon_S e^{-\varepsilon_S \gamma_{jh} }$.
We conclude that the GML approach is equivalent to Algorithm \ref{algo:RWS} wherein Step 6 is now replaced by (\ref{rew_ML}). Finally, it is not difficult to see that the sequence $\hat X^{(l)}$ generated by the GML method converges to the set of stationary points of the Problem (\ref{pb_MM}) where the log-sum penalty now is
\al{f(X)&:=(2n+1)\sum_{j>h}\log ( q_{jh}(\Dd(X))+\varepsilon_S )\nn\\ &+(n+1)\sum_{j=1}^m\log ( q_{jj}(\Dd(X))+\varepsilon_S ).\nn}

\section{Identification of Latent-variable Graphical Models} \label{sec:Bayes_SL}
Consider a Gaussian discrete-time zero mean full rank stationary stochastic process $z=\{\,z(t),\; t\in\mathbb Z\,\}$. $z$ is composed by  $m$ manifest variables and $r$ latent variables, so that $z=[\, (y)^T\,(x)^T\,]^T$ with $y=[\,y_1 \ldots y_m\,]^T$ and $x=[\,x_{1} \ldots x_{r}\,]^T$. We assume that $y$ is an AR process of order $n$ (known) and dimension $m$ (known). $r$ as well as the PSD of $z$ are unknown. We assume to collect the data $ y^N:=\{\,  y(1), y(2)\ldots$ $  y(N)\,\}$ from the manifest process $y$. We want to estimate the PSD of $y$, say $\Phi$, in such a way that it corresponds to a latent-variable graphical model. The latter is an interaction graph with two layers: latent nodes $x_1\ldots x_r$ are in the upper level, manifest nodes $y_1\ldots y_m $ are in the lower level, $r\ll m$ and there are few edges among the manifest nodes. An example of a latent-variable graphical model is provided in Figure \ref{fig:graph_ex} (right). The powerfulness of latent-variable graphical models is that the introduction of few latent variables (in respect to the manifest ones) may reduces drastically the conditional dependences among the manifest variables. Accordingly, in such graphs we expect that the interdependence relations among the manifest variables are mainly explained by few and common latent variables. In \cite{LATENTG} it has been shown that $\Phi$ corresponds to a latent-variable graphical model if it admits the following decomposition
\al{\Phi^{-1}=\Sigma-\Lambda}
where \al{\label{def_Sigma_Lambb}\Sigma&=S_0+\frac{1}{2}\sum_{k=1}^n  S_k e^{-i\theta k} +   S_k^Te^{i\theta k}\nn\\ \Lambda&=L_0+\frac{1}{2}\sum_{k=1}^n  L_k e^{-i\theta k} +   L_k^Te^{i\theta k}.} 
$\Sigma\succ 0$ is sparse and its support reflects the conditional dependence relations among the manifest variables. $\Lambda\succeq 0$ is low rank and its rank is equal to $r$, i.e the number of latent variables. The approximate negative log-likelihood of
  $  y(n
+1) \ldots  y(N)$ given $  y(1) \ldots   y(n)$ under the aforementioned model is: 
\al{\label{cond_PDF_SL}\ell(   y^N; \Sigma,\Lambda):=\ell(   y^N; \Sigma-\Lambda)} 
where $\ell(   y^N;\cdot)$ has been defined as in (\ref{cond_PDF}). Accordingly, the regularized ML  estimator of $\Sigma$ and $\Lambda$ is  \citep{LATENTG}:
\al{\label{pb_SL_2016}(\hat \Sigma,\hat \Lambda)= &\underset{\Sigma, \Lambda \in \Qc_{m,n}}{ \mathrm{argmin}} \ell( y^N; \Sigma,\Lambda)+ \gamma_S\bar h_{\infty}(\Sigma)+\gamma_L \bar r_\star(\Lambda) \nn\\
& \hbox{ s.t. }  \Sigma\succ0 , \; \; \Lambda \succeq 0.}
Here, $\bar r_\star(\Lambda)=\tr\int \Lambda$ is the nuclear norm of $\Lambda \in \Qc_{m,n}$ with $\Lambda\succeq 0$.  $\gamma_S,\gamma_L>0$ are the regularization parameters for $\Sigma$ and $\Lambda$, respectively. The term $\bar h_\infty(\Sigma)$, defined in (\ref{h_inf}), induces sparsity on $\Sigma$. 
The term $\bar r_\star(\Lambda)$ induces low rank on $\Lambda$, indeed it has been shown that such a function is the convex envelop of $\mathrm{rank} (\Lambda)$, see Proposition 3.1 in \cite{LATENTG}. It is worth noting that the estimate of $r$ is given by the numerical rank of $\hat \Lambda$. Since $\Sigma-\Lambda$ and $\Lambda$ belong to $\Qc_{m,n}$ and in view of (\ref{set_Qc_mn_reparametrized}), we can define $X,H\in\Qb_{m(n+1)}$ such that
 \al{\label{XL_par}\Sigma-\Lambda=\Delta X \Delta^*, \; \; \Lambda=\Delta H \Delta^*}
 and $X,H\succeq 0$. It is possible to prove that Problem (\ref{pb_SL_2016}) can be rewritten in terms of $X$ and $H$
as follows: 
\al{\label{pb_SL_2016_ref}(\hat X,\hat H)= &\underset{X, H \in \Qb_{m(n+1)}}{ \mathrm{argmin}} \ell( y^N; X)+ \gamma_S h_{\infty}(\Dd(X+H))\nn\\
& \hspace{4.4cm}+\gamma_L \tr(H)\nn\\
& \hbox{ s.t. }  X_{00}\succ 0,\; \; X\succeq 0 , \; \; H\succeq 0}
where $\ell( y^N; X)$ does coincide with the first term in the objective function of (\ref{pb_td2}). Also in this case the optimal solution $(\hat \Sigma,\hat\Lambda)$
depends on the regularization parameters $\gamma_S$, $\gamma_L$.   The values
$(\gamma_S,\gamma_L)$ can be selected by considering a 2-dimensional grid. Then, the candidate models can be ranked by using a BIC criterium similar to the one in \cite{SONGSIRI_TOP_SEL_2010}. The latter measures the complexity by thresholding the partial coherence of the sparse component and the singular values of the low rank component. Alternatively,  \cite{LATENTG} proposes a score function based on the Kullback-Leibler divergence to rank the candidate models.

Similarly to the case without latent variables, the weaknesses of the regularization Problem (\ref{pb_SL_2016}) are that: (i) the nonnull entries in $\Sigma$ are penalized in a different way; (ii) the nonnull eigenvalues of $\Lambda$, understood as functions over the unit circle, are penalized in a different way.  
The sparse regularizer can be made ``democratic'' in the same way of before. Regarding the low rank part,  we can replace $\bar r_\star(\Lambda)$ with the ``weighted'' penalty:
\al{ \bar r_W (\Lambda)=\tr \left(Q\int \Lambda\right)\nn}
where the weight matrix $Q\in\Qb_m$ is such that $Q\succ 0$. Thus, we consider the problem:
\al{\label{pb_SL_rw}(\hat \Sigma,\hat \Lambda)= &\underset{\Sigma, \Lambda \in \Qc_{m,n}}{ \mathrm{argmin}} \ell( y^N; \Sigma,\Lambda)+ \bar h_{W}(\Sigma)+ \bar r_W(\Lambda) \nn\\
& \hbox{ s.t. }  \Sigma\succ0 , \; \; \Lambda \succeq 0.} Using the parametrization in (\ref{XL_par}), 
we have that
\al{r_W(\Lambda)&=\tr\left(Q\int \Delta H \Delta^*\right)=\tr\left(\int \Delta^* Q\Delta H\right)\nn\\
&=\tr((I_{n+1}\otimes Q )H)}
where we exploited the well known identity $\int e^{i k \theta}=\delta_k$ and $\delta_k$ denotes the Kronecker delta function. Thus, we consider the problem:
\al{\label{pb_SL_rw2}(\hat X,\hat H)= &\underset{X,H\in \Qb_{m(n+1)}}{ \mathrm{argmin}} \ell ( y^N; X)+ h_{W}(\Dd(X+H))\nn\\
&\hspace{4.2cm}+\tr((I\otimes Q)H)\nn\\ & \hspace{0.4cm}\hbox{ s.t. }  X_{00}\succ 0, \; \;X\succeq 0,\; \; H\succeq 0.} 
Note that, (\ref{pb_SL_rw}) and (\ref{pb_SL_rw2}) are equivalent provided that $\Delta \hat X\Delta^*\succ 0$. 
\teo
\label{teo_SL_exist}The dual problem of (\ref{pb_SL_rw2}) is  
\al{\label{dual_SLrw2} &\underset{\substack{Z\in \Mb_{m,n}\\ W\in \Qb_{m}}}{ \mathrm{max}} \log | W|+m\nn\\ & \hbox{ s.t. } W\succ 0\nn\\& \hspace*{0.7cm}\Td(\hat R)+\Td(Z)\succeq \left[\begin{array}{cc}W &  0 \\0  & 
0\end{array}\right]\nn\\ & \hspace*{0.7cm} \sum_{k=0}^n |(Z_k)_{jh}|+|(Z_k)_{hj}|\leq \frac{2\gamma_{jh}}{N-n},\; \; j>h\nn\\
 &\hspace*{0.7cm} \sum_{k=0}^n |(Z_k)_{jj}|\leq \frac{2\gamma_{jj}}{N-n},\; \; j=h\nn\\
 & \hspace*{0.7cm} I\otimes Q +\Td(Z)\succeq 0.} The latter does admit solution $(\hat X,\hat H)$ such that $\hat X$ is unique and $\Delta \hat X \Delta^*\succ 0$. Accordingly (\ref{pb_SL_rw}) and (\ref{pb_SL_rw2}) are equivalent.\eteo

Theorem \ref{teo_SL_exist} is important not only for establishing the existence of the solution of (\ref{pb_SL_rw}) and (\ref{pb_SL_rw2}), but it provides also a tractable formulation for the computation of the solution, see Section \ref{sec:imple}. It is worth noting that Problem (\ref{pb_SL_rw}) may have more than one solution. On the other hand, if we compute an optimal solution of (\ref{dual_SLrw2}), then from the latter we can recover the solution of (\ref{pb_SL_rw2}) by solving a system of linear equations in $H$ (a similar idea has been used in Section III.C in \cite{LATENTG}). The uniqueness of the solution of this system of linear equations is guaranteed provided that: (i) $Q$ has a sufficient number of eigenvalues which are sufficiently large; (ii) there is a sufficient number of $\gamma_{jh}$s with $j\neq h$ taking sufficiently large values. The latter implies the uniqueness of the solution to (\ref{pb_SL_rw2}) and thus the uniqueness of the one to (\ref{pb_SL_rw}).

To select a suitable set of $\gamma_{jh}$ and a suitable $Q$, we design a reweighted algorithm which exploits a Bayesian perspective. We model $\Sigma$ and $\Lambda$ as stochastic processes taking values in $\Qc_{m,n}^+$ and in the closure of $\Qc_{m,n}^+$, respectively. Let $p(\Sigma)$ and $p(\Lambda)$ denote the prior of $\Sigma$ and $\Lambda$, respectively. We assume that $\Sigma$ and $\Lambda$ are independent, that is the joint PDF of $\Sigma$ and $\Lambda$ is such that $p(\Sigma,\Lambda)=p(\Sigma)p(\Lambda)$. We set the prior  of $\Sigma$ as in (\ref{p_Sigma}). Regarding $\Lambda$, $L_k$s are modeled as independent random matrices. More precisely, we model $L_0$ as a
Wishart random matrix with $m+1$ degrees of freedom and variance $mQ^{-1}$:
\al{p_Q( L_0)=\frac{e^{- \tr\left(Q   L_0\right)}}{\sqrt{ 2^{(m+1)m}|Q|^{-(m+1) }}G_m(\frac{m+1}{2})}}
where $Q\in \Qb_m$ such that $Q\succ 0$ and $G_m$ is the multivariate gamma function. Finally, we attach an uninformative prior on $L_k$ with $k\geq 1$. 
Then, the negative log-likelihood of $y^N$, $\Sigma$ and $\Lambda$ is defined as
\al{\ell(& y^N, \Sigma,\Lambda; \pmb \gamma, Q)= -\log p_{\pmb \gamma,Q}( y^N,\Sigma,\Lambda)\nn \\
&= -\log p( y^N|\Sigma,\Lambda)-\log p_{\pmb \gamma}(\Sigma)-\log p_{Q}(\Lambda)+c\nn} where $c$ is a constant term not depending on $y^N$, $\Sigma$, $\Lambda$, $\pmb{\gamma}$ and $Q$. Note that, $-\log p( y^N|\Sigma,\Lambda)$ does coincide with $\ell(   y^N; \Sigma ,\Lambda)$ defined in (\ref{cond_PDF_SL}). Also in this case, it is not possible to find an analytical expression for the negative log-marginal likelihood. Therefore, we have the following upper bound for $\ell( y^N , \Sigma,\Lambda; \pmb \gamma, Q)$:
\al{\tilde \ell( y^N &, \Sigma,\Lambda; \pmb \gamma, Q)=
\ell(   y^N; \Sigma ,\Lambda)+\bar h_W(\Sigma)+\sum_{j\geq h} \log c_{jh}\nn\\ &+\bar r_W(\Lambda)-\frac{m+1}{2}\log|Q|+c.} Then, according to the GML approach we consider the following two-step algorithm:
\al{\label{2step_SL}(\hat \Sigma^{(l)},\hat \Lambda^{(l)}) =& \underset{  \Sigma,\Lambda}{\;\mathrm{argmin}}\;  \tilde \ell ( y^N,\Sigma,\Lambda;\hat{\pmb \gamma}^{(l)},\hat Q^{(l)})\\ 
\label{2step_gammaQ}(\hat{\pmb \gamma}^{(l+1)},\hat Q^{(l+1) })=& \underset{\pmb\gamma,Q }{\;\mathrm{argmin}}\;   \tilde \ell ( y^N,\hat \Sigma^{(l)},\hat \Lambda^{(l)};\pmb \gamma, Q).} 
Clearly, (\ref{2step_SL}) is the MAP estimator of $\Sigma$ and $\Lambda$
given the current choice of $\pmb \gamma$ and $Q$, while (\ref{2step_gammaQ}) is the estimator of  $\pmb \gamma$ and $Q$ using the current MAP estimate of $\Sigma$ and $\Lambda$ as a direct observation. The objective function in (\ref{2step_gammaQ}) can be split in two terms depending on $\pmb \gamma$ and $Q$, respectively. Accordingly, $\hat{\pmb \gamma}^{(l+1)}$
takes a form similar to the one in (\ref{rew_MLsenza_e}). Regarding $Q$, we have:
\al{\label{opt_Q}\hat Q^{(l+1)}=& \underset{Q\succ 0 }{\;\mathrm{argmin}}\;   \tr(Q\hat L_0^{(l)})-\frac{m+1}{2}\log|Q|.}
\pro \label{prop_opt_Q}Under the assumption that $\hat L_0^{(l)}\succ 0$, we have that
\al{\label{opt_Q2}\hat Q^{(l+1)}=\frac{m+1}{2}(\hat L_0^{(l)})^{-1}.}
\epro

To deal also with the case $\hat L_0^{(l)}$ singular matrix we consider the modified updating rule:
\al{\label{opt_Q3}\hat Q^{(l+1)}=\frac{m+1}{2}(\hat L_0^{(l)}+\varepsilon_L I)^{-1}}
where $\varepsilon_L>0$. It is not difficult to see that the latter modification is equivalent to assume that $Q$ is modeled as a Wishart random matrix  with $m+1$ degrees of freedom and variance $m\varepsilon_L^{-1}I$, i.e. $p(Q)=e^{-\varepsilon_L\tr(Q)}/(\sqrt{2^{(m+1)m}\varepsilon_L^{-m(m+1)}}G_m((m+1)/2))$, which is assumed to be independent of $\pmb \gamma$. Algorithm \ref{algo:RWSL} describes the corresponding procedure which is clearly an iterative reweighted algorithm. Here, we recall that $\hat S^{(l)}=\Dd(\hat X^{(l)}+\hat H^{(l)})$ and $\hat L_0^{(l)}=\Dd_0(\hat H^{(l)})$.

\begin{algorithm}
\caption{Iterative Reweighted Algorithm}
\label{algo:RWSL}
\begin{algorithmic}[1] \small
\STATE $l=0$
\STATE Initialize $\hat \gamma_{jh}^{(0)}$ with $j,h=1\ldots m$ and $\hat Q^{(0)}$
\REPEAT 
\STATE Solve Problem (\ref{pb_SL_rw2}) with $\pmb{\gamma}=\hat{\pmb{\gamma}}^{(l)}$ and $Q=\hat Q^{(l)}$; let  \STATEx \hspace{0.3cm}$(\hat X^{(l)},\hat H^{(l)})$ denote the corresponding solution
\STATE Update the weights 
\al{\hat{\gamma}^{(l+1)}_{jh}&=\left\{\begin{array}{ll} \frac{n+1}{q_{jh}(\Dd(\hat X^{(l)}+\hat H^{(l)}))+\varepsilon_S },& \hbox{ if } j=h  \\ \frac{2n+1}{q_{jh}(\Dd(\hat X^{(l)}+\hat H^{(l)}))+\varepsilon_S  }, & \hbox{ if } j>h \end{array}\right. \nn \\
\hat Q^{(l+1)}&=\frac{m+1}{2}(\Dd_0(\hat H^{(l)})+\varepsilon_L I)^{-1}\nn}
\STATE $l\leftarrow l+1$
\UNTIL{ $\|\hat X^{(l)}-\hat X^{(l-1)}\|+\|\hat H^{(l)}-\hat H^{(l-1)}\|\geq \varepsilon$ \textbf{or} $l\leq l_{\mathrm{MAX}}$ }
\STATE $\hat \Sigma=\Delta(\hat X^{(l)}+\hat H^{(l)})\Delta^*$, $\hat \Lambda=\Delta\hat H^{(l)}\Delta^*$
\STATE  $\hat \Phi=(\hat \Sigma-\hat \Lambda)^{-1}$
\end{algorithmic}
\end{algorithm}

It is not difficult to see that Algorithm \ref{algo:RWSL} can be interpreted as the MM algorithm 
for solving the regularized ML problem:
\al{\label{pb_logdet_SL}& \underset{\substack{X,H\in \Qb_{m(n+1)}}}{ \mathrm{argmin}} \; \; f(X,H)
\nn\\ & \hbox{ s.t. }  X_{00}\succ 0,\; \; X\succeq 0, \; \; H\succeq 0}
where
\al{f(X,H)&:=(2n+1)\sum_{j>h}\log ( q_{jh}(\Dd(X+H))+\varepsilon_S )\nn\\ &+(n+1)\sum_{j=1}^m\log ( q_{jj}(\Dd(X+H))+\varepsilon_S )\nn\\
&+\frac{m+1}{2}\log |\Dd_0(H)+\varepsilon_L I |+\ell ( y^N; X).\nn}
Beside the sparse-inducing part (already analyzed in Section \ref{sec:Bayes_S}), it is known that the log-det penalty outperforms the nuclear norm for estimating the correct rank, see \cite{fazel2003log}. 

\corr \label{cor_conv_RWSL}Consider the sequence $(\hat X^{(l)},\hat H^{(l)})$ generated by Algorithm \ref{algo:RWSL} and assume that the solution of (\ref{pb_SL_rw2}), with weights $\hat{\pmb \gamma}^{(l)}$ and $\hat Q^{(l)}$, is always unique. Then, $(\hat X^{(l)},\hat H^{(l)})$ converges to the set of stationary points of Problem (\ref{pb_logdet_SL}).
\ecorr

\begin{figure}[htbp]
\centering
\subfloat{ \includegraphics[width=0.98\columnwidth]{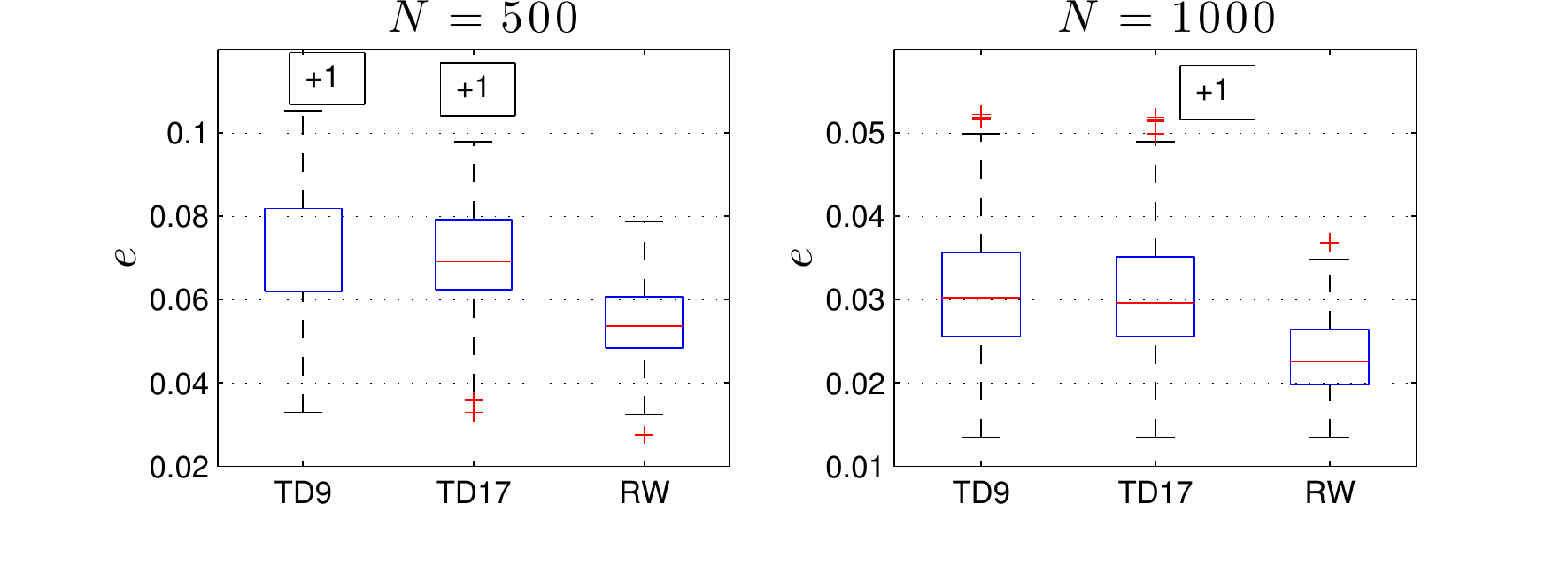}} \;
\subfloat{ \includegraphics[width=0.98\columnwidth]{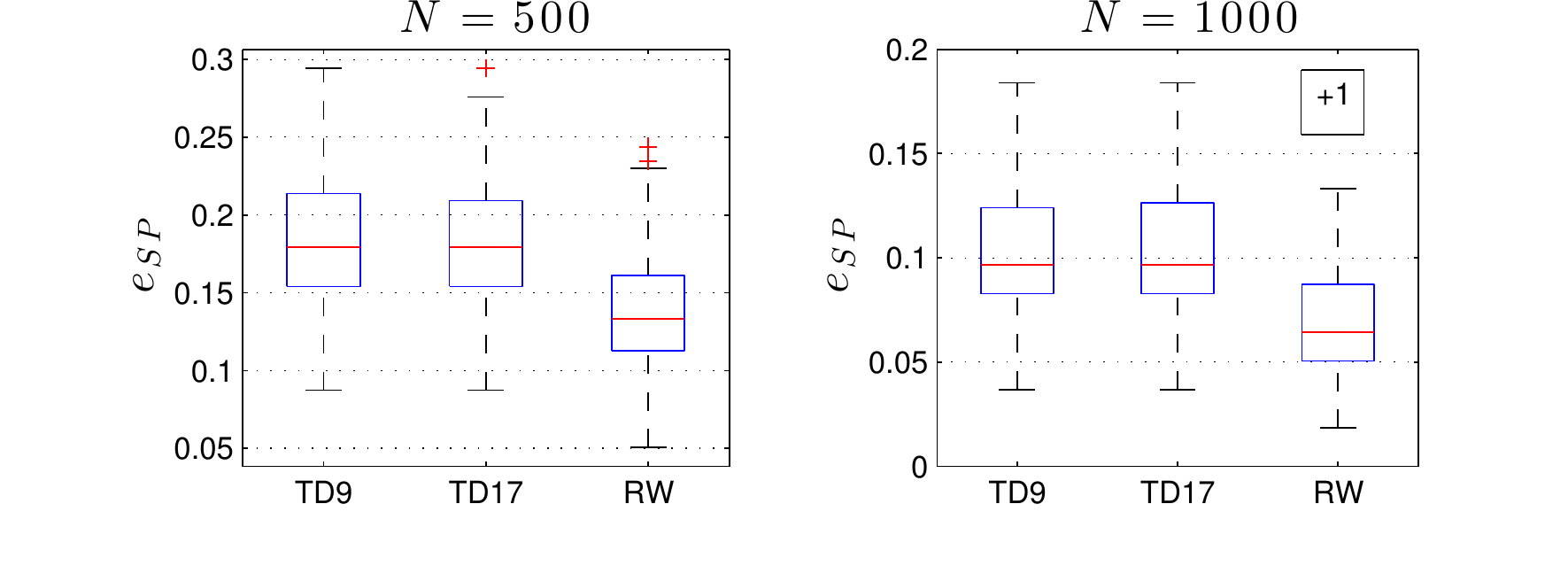}} 
\caption{First experiment with sparse AR models of order $n=1$. 
} \label{Fig:exp1}
\end{figure}
\begin{figure}[htbp]
\centering
\subfloat{ \includegraphics[width=0.98\columnwidth]{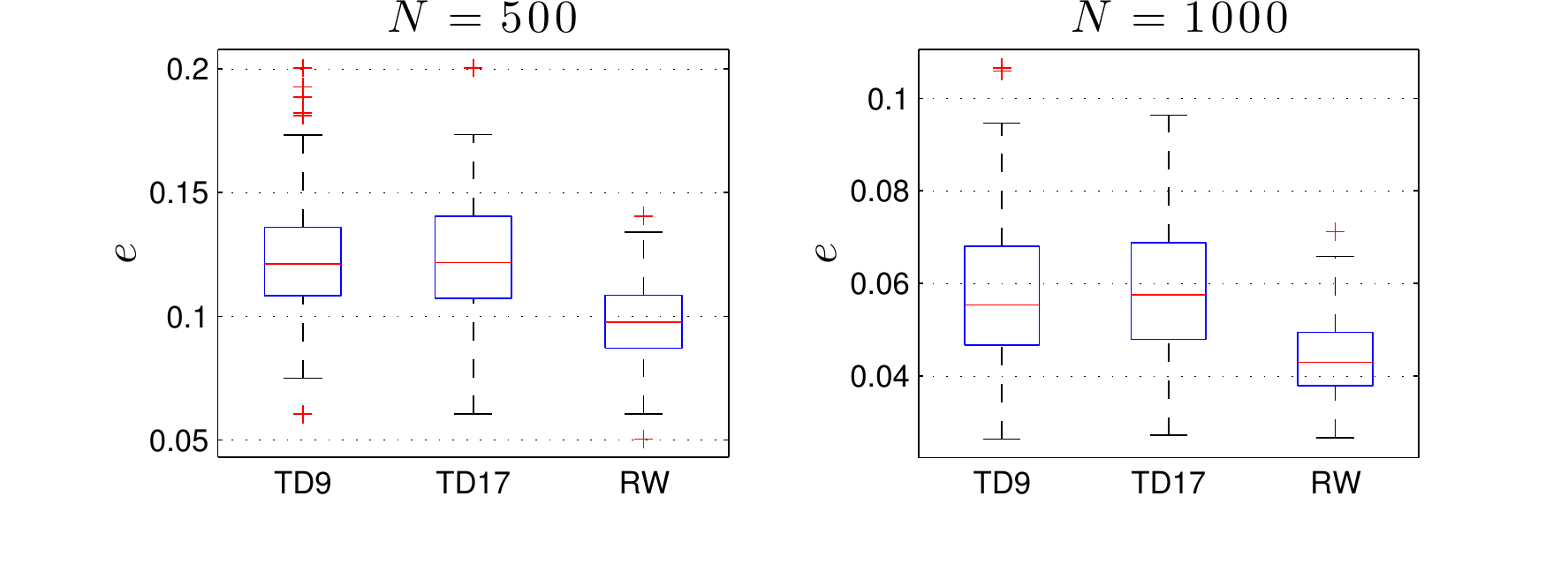}} \;
\subfloat{ \includegraphics[width=0.98\columnwidth]{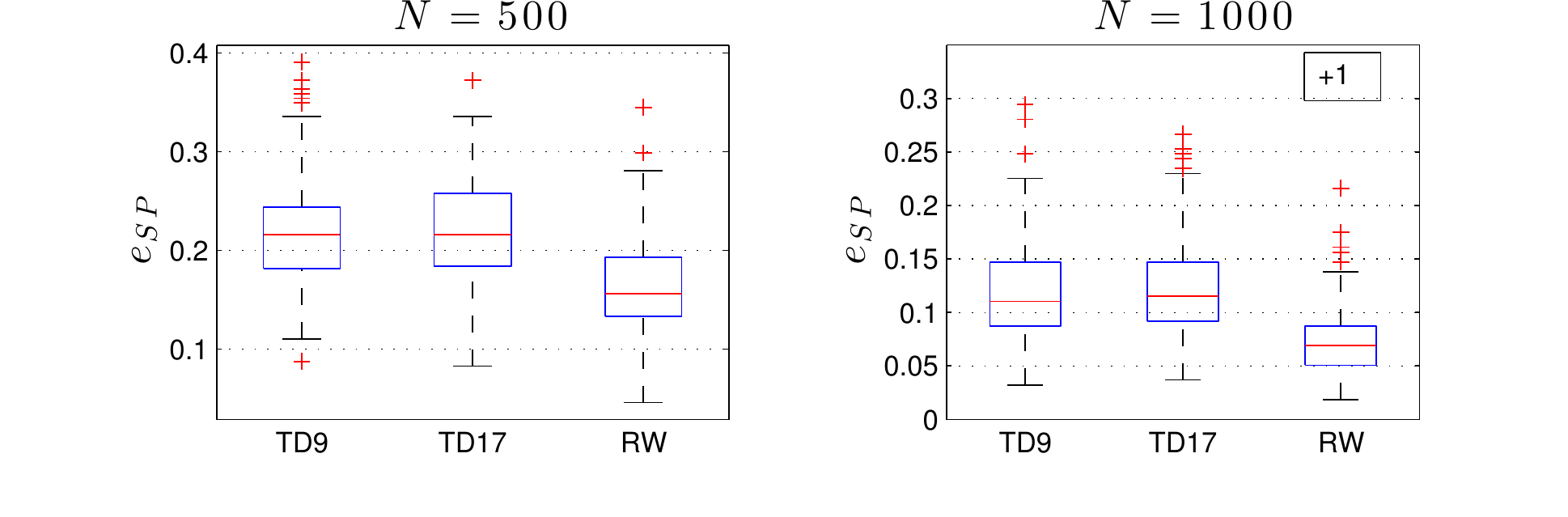}}
\caption{Second experiment with sparse AR models of order $n=2$. 
} \label{Fig:exp2}
\end{figure}
\begin{figure}[htbp]
\centering
\subfloat{ \includegraphics[width=0.98\columnwidth]{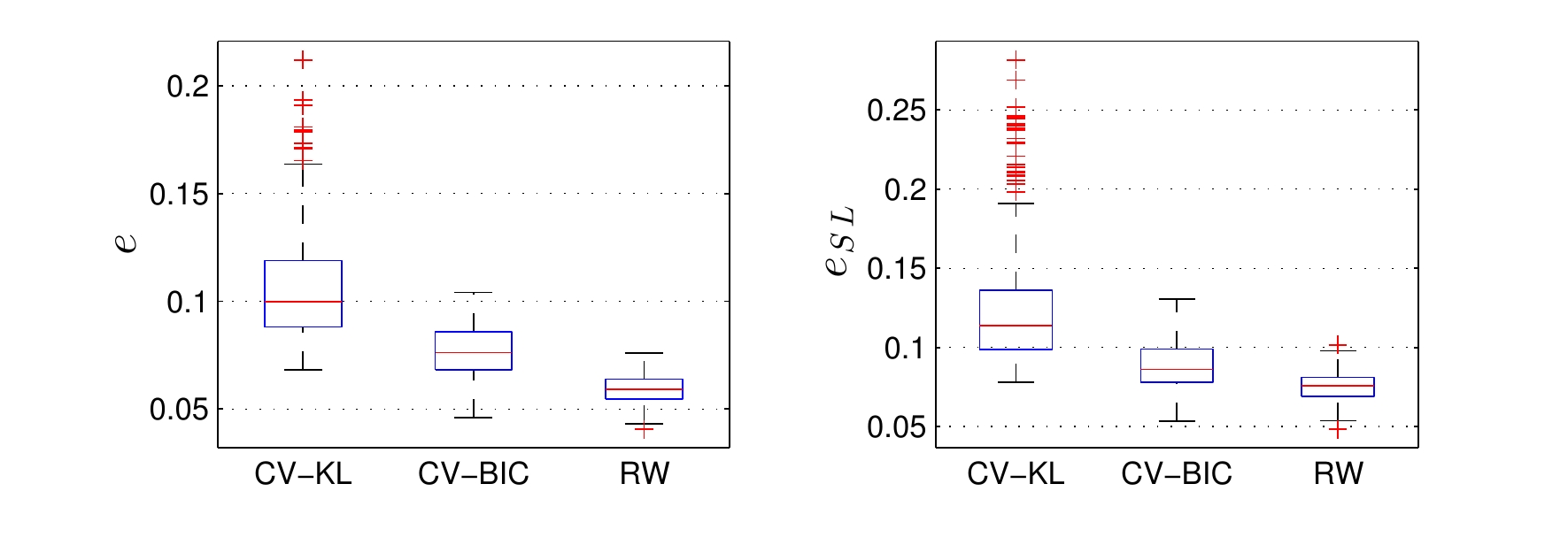}} 
\caption{First experiment with latent-variable AR models of order $n=2$ with $s=0.1$ and $r=2$. 
} \label{Fig:exp3}
\end{figure}
\begin{figure}[htbp]
\centering
\subfloat{\includegraphics[width=0.98\columnwidth]{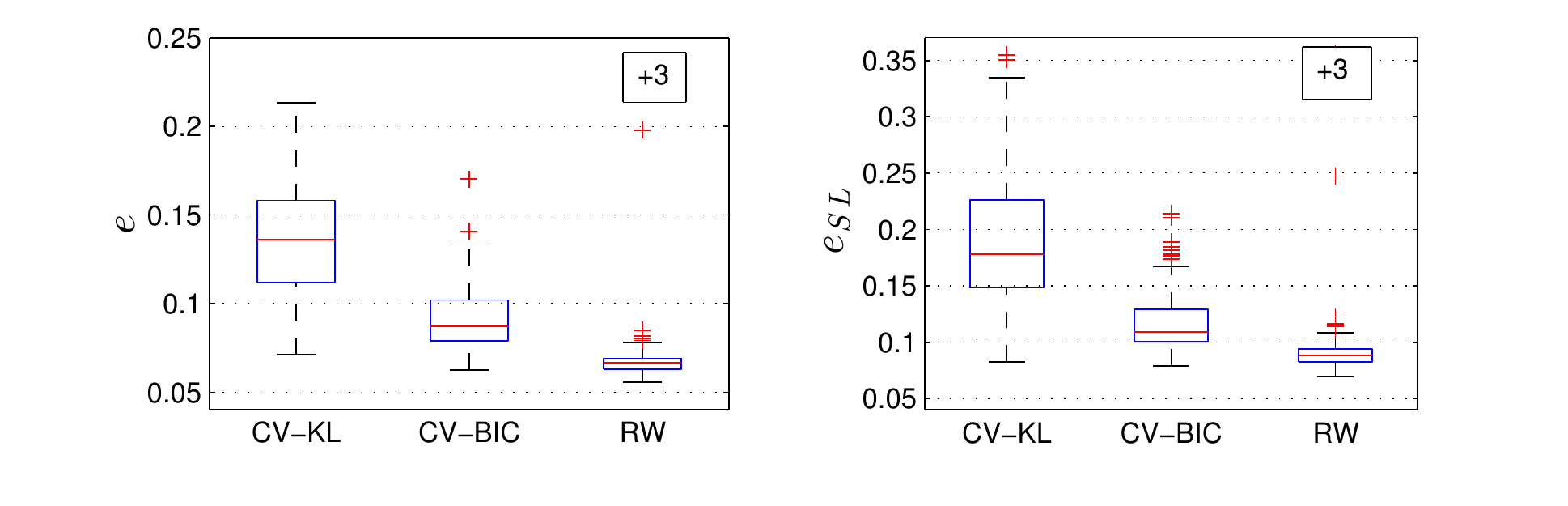}} 
\caption{Second experiment with latent-variable  AR models of order $n=2$ with $s=0.1$ and $r=5$. 
 } \label{Fig:exp4}
\end{figure}

\section{Simulation Results}\label{sec:sim}
In this section we test the performance of the proposed  empirical Bayes estimators for learning AR graphical models. 

\subsection{Implementation details}\label{sec:imple}

The empirical Bayes estimate of a sparse AR graphical model is obtained by running Algorithm \ref{algo:RWS} wherein Step 6 is replaced by the updating rule (\ref{rew_ML}). The most challenging part is the computation of $\hat X^{(l)}$ (Step 4) which is obtained through the dual problem (\ref{dual_rw2}). The optimal solution of the latter is computed by implementing a first-order projection algorithm. The empirical  Bayes estimate of a latent-variable AR graphical model is obtained by running Algorithm \ref{algo:RWSL}. Here, the most challenging part 
is the computation of $(\hat X^{(l)},\hat H^{(l)})$, in Step 4, which is obtained through the dual problem (\ref{dual_SLrw2}). The optimal solution of the latter is computed by using the alternating direction method of multipliers (ADMM). More precisely, we introduce an auxiliary variable $K=I\otimes Q+\Td(Z)$. In the first step, we minimize the augmented Lagrangian with respect to $W$ and $Z$ under the first four constraints of (\ref{dual_SLrw2}). This minimization is performed by using a first-order projection algorithm. In the second step, we minimize the augmented Lagrangian with respect to $K$ under the last constraint of  (\ref{dual_SLrw2}). This optimal solution can be computed in closed form. Finally, the corresponding Matlab functions are available at \url{https://github.com/MattiaZ85/EBL-AR-GM}.

\subsection{Identification of sparse graphical models}
We consider two Monte Carlo experiments which are structured as follows.
We generate 200 AR models of dimension $m=30$, of order $n$ and  whose fraction of nonnull entries in the inverse PSD is equal to $0.1$. The position of such nonnull entries is chosen randomly for each model. For each model we generate a finite data sequence of length $N$  and we consider the following estimators:
\begin{itemize}
\item TD9 is the estimator proposed in \cite{SONGSIRI_TOP_SEL_2010} where the number of points for tracing the trade-off curve is set equal to $9$ and the threshold for the partial coherence is set equal to $0.1$;
\item TD17 is the estimator proposed in \cite{SONGSIRI_TOP_SEL_2010} where the number of points for tracing the trade-off curve is set equal to $17$ and the threshold for the partial coherence is set equal to $0.1$;
\item RW  is the empirical Bayes estimator that we have proposed in Section \ref{sec:Bayes_S} where $\varepsilon_S=10^{-3}$. The hyperparameters vector $\pmb{\gamma}$ is initialized with (\ref{rew_ML})
where $\hat S^{(l)}=\Dd(X_B)$, $X_B\in\Qb_{m(n+1)}$ is such that $\hat \Phi_B=(\Delta X_B\Delta^*)^{-1}$ and $\hat \Phi_B$ is the Burg estimator computed from the data\footnote{It is worth recalling that the Burg method provides an estimate of the PSD which corresponds to an AR model of a certain order, fixed by the user. In our case, the latter is fixed equal to $n$; in this way $\hat \Phi_B=(\Delta X_B \Delta^*)^{-1}$.}. The threshold for the stopping condition is set equal to $\varepsilon=10^{-4}$ and the maximum number of iterations is set equal to $l_{\mathrm{MAX}}=50$.\end{itemize}
For each estimator:
\begin{itemize}
\item we compute the relative error in the estimation of the inverse PSD
coefficients as $e=\|\hat S-S\|^2 /\|S\|^2$
where $\hat S=[\,\hat S_0\,\hat S_1\,\ldots \,\hat S_n]\in \Mb_{m,n}$ are the coefficients of the inverse of $\hat \Phi$ (estimated PSD), while $S=[\, S_0\, S_1\,\ldots \, S_n]\in\Mb_{m,n}$ are the ones of the inverse of $\Phi$ (true PSD); 
\item  we compute the fraction of null and nonnull misplaced entries in the inverse of the estimated PSD with respect to the inverse of the true PSD as $e_{SP}= \#\hat E_m /(m(m-1)/2)$
 where $\#\hat E_m$ denotes the total number of null and nonnull misplaced entries in the inverse of the estimated PSD.   
\end{itemize}

In the first experiment, the AR models are of order $n=1$ and we have considered  two different lengths of the data, that is $N=500$ and $N=1000$. The boxplots of $e$ and $e_{SP}$ are depicted in Figure \ref{Fig:exp1}. As we can see, RW outperforms TD9 and TD17 both in terms of estimation error and of misplaces entries. It is worth noting that TD9 and TD17 perform in the same way: this confirms the observation that for tracing accurately the trade-off curve it is required just a small number of points \citep{SONGSIRI_TOP_SEL_2010}. We have noticed that TD9 and TD17 estimate the inverse PSD too sparse in respect to the true one. For this reason, we also have considered TD9 and TD17 with the threshold for the partial coherence equal to $0.05$: the resulting performances are worse than the ones with $0.1$. In the second experiment, the AR models are of order $n=2$ and we have considered the two different lengths of the data as before. The boxplots of $e$ and $e_{SP}$ are depicted in Figure \ref{Fig:exp2}. The results are similar to the ones of the previous experiment: RW outperforms the other two estimators; TD9 and TD17 perform in the same way.

\subsection{Identification of latent-variable AR graphical models}

We consider two Monte Carlo experiments which are structured as follows.
We generate 200 AR models of dimension $m=30$ of order $n=2$ whose inverse PSD has a sparse plus low rank decomposition, i.e. $\Phi^{-1}=\Sigma-\Lambda$ with $\Sigma$ sparse and $\Lambda$ low rank. The fraction of nonnull entries in the sparse part is denoted by $s$; the position of the nonnull entries is chosen randomly for each model. The rank of $\Lambda$ is denoted by $r$; for each model, the image of $\Lambda$ is chosen randomly in such a way that $\Sigma-\Lambda$ is positive definite over the unit circle. For each model we generate a finite data sequence of length $N=1000$  and we consider the following estimators:
\begin{itemize}
\item TD-KL is the estimator proposed in \cite{LATENTG} where the model is chosen using the score function based on the Kullback-Leibler divergence; the threshold for the partial coherence (sparse component) is set equal to $0.1$; the threshold for the normalized singular values (low rank component) is set equal to $0.1$; the regularization grid (16 points) is chosen in such a way that the candidate model, say $\hat \Phi=(\hat S-\hat L)^{-1}$, ranges from $\hat S$ diagonal - $\mathrm{rank}(\hat L)$ large to $\hat S$ moderately sparse and $\mathrm{rank}(\hat L)=0$.

\item TD-BIC is the estimator proposed in \cite{LATENTG} where the model is chosen using the BIC score function; the threshold parameters and the regularization grid are chosen as in the previous estimator. 
\item RW  is the empirical Bayes estimator that we have proposed in Section \ref{sec:Bayes_SL} where $\varepsilon_S=10^{-3}$ and $\varepsilon_L=10^{-3}$. The hyperparameters are initialized as follows. Let $\hat \Phi_B=(\Delta X_B \Delta^*)^{-1}$ be the Burg estimator of the PSD computed from the data. Then, the initial value for $\pmb{\gamma}$ is given by (\ref{rew_ML})
where $\hat S^{(l)}=(1+\alpha) \Dd(X_B)$, the initial value for $Q$ is given by (\ref{opt_Q3}) where $\hat L^{(l)}=\alpha\Dd_0(X_B)$ and $\alpha=0.1$. The threshold for the stopping condition is set equal to $\varepsilon=10^{-4}$ and the maximum number of iterations is set equal to $l_{\mathrm{MAX}}=50$.\end{itemize} 

For each estimator:
\begin{itemize}
\item we compute the relative error in the estimation of the inverse PSD
coefficients as $e=\|\hat S-\hat L-(S-L)\|^2 /\|S-L\|^2$ 
where $\hat S=[\, \hat S_0\,\hat S_1\,\ldots \,\hat S_n]\in\Mb_{m,n}$ are the coefficients of $\hat \Sigma$, while $\hat L=[\,\hat L_0\,\hat L_1\,\ldots \,\hat L_n]\in\Mb_{m,n}$ are the ones of $\hat\Lambda$. $\hat \Phi=(\hat \Sigma-\hat \Lambda)^{-1}$ denotes the estimated PSD. In a similar way, $S$ and $L$ contain the coefficients of $\Sigma$ and $\Lambda$ such that $\Phi=(\Sigma-\Lambda)^{-1}$ is the true PSD.
\item we compute the relative error in the estimation of the sparse and low rank coefficients in the inverse PSD as $e_{SL}=(\|\hat S-S\|^2+\|\hat L-L\|^2) /\|S-L\|^2$;
\item we compute the complexity of the estimated model as  
$C=\#p /(m(m+1)/2+m^2 n)$ where $\#p$ is the number of parameters needed to characterize the model $\hat \Phi=(\hat \Sigma-\hat \Lambda)^{-1}$. For instance if $\hat \Sigma$ has $s$ nonnull entries and $\mathrm{rank}(\Lambda)=r$, then $\#p=(s+m)/2+sn+rm(n+1)$. The denominator in $C$ is the number of coefficients needed in an unstructured model, i.e. the corresponding interaction graph is complete.
\end{itemize}

\begin{figure}
\begin{center}
\begin{tabular}{|l|c|}
\hline
 1st exp.& $\bar C$ \\ 
\hline              
TRUE & 0.1890 \\
 TD-KL  &  0.0902  \\
  TD-BIC  & 0.1012    \\ 
   RW  &  0.2100 \\
\hline
\end{tabular}
\quad\quad
\begin{tabular}[1.0\textwidth]{|l| c|}
\hline
2nd exp. & $\bar C$  
\\  
\hline               
TRUE & 0.2684 \\
 TD-KL  &   0.1029 \\
  TD-BIC  &  0.1158   \\
   RW  & 0.2810 \\
   \hline
\end{tabular}
\bigskip\caption{Average complexity of the estimated models in the first Monte Carlo experiment where $s=0.1$, $r=2$ (left table) and in the second Monte Carlo experiment where $s=0.1$, $r=5$ (right table).} 
\label{tab:rank}
\end{center}
\end{figure} 

In the first experiment, the AR models are with $s=0.1$ and $r=2$. The relative errors $e$ and $e_{SL}$ are depicted in Figure \ref{Fig:exp3}. As we can see, RW provides the best performance while TD-KL the worst performance. The average complexity of the estimated models for each estimator (and denoted by $\bar C$) is shown in 
the table on the left of Figure \ref{tab:rank}. As we can see, the average complexity of the models estimated by RW is closer to the true one than the one of TD-KL and TD-BIC. More precisely, TD-KL and TD-BIC estimate models which are extremely simple, in respect to the true complexity, and the price to pay is the inferiority in respect to $e$ and $e_{SL}$. In the second experiment, the AR models are with $s=0.1$ and $r=5$. The relative errors $e$ and $e_{SL}$ are depicted in Figure \ref{Fig:exp4}, while the average complexity of the estimated models for each estimator is shown in the table on the right of Figure \ref{tab:rank}. As we can see, the conclusions of the previous case still hold. Finally, we have tested RW for different values of $\alpha$ such that $0.08\leq \alpha\leq 0.15$. We have noticed that RW provides similar performances of before.

\section{Conclusion} \label{sec:concl}
We have analyzed the problem of estimating sparse and latent-variable AR graphical models. These two problems are traditionally solved by using $\ell_1$-like and nuclear norm-like regularizes. The latter introduce a magnitude imbalance in the optimization problem which produces an estimator with a remarkable mean squared error. We have proposed two empirical Bayes estimators which counteract such an imbalance. The hyperparameters of these estimators are computed by the generalized maximum likelihood approach which leads to an iterative reweighted method. Simulation results showed the benefit to introduce a prior for the model that we have to estimate.

\section*{Appendix} 
 
 \subsection*{Proof of Theorem \ref{teo_S1}}
 First, we can rewrite (\ref{pb_rw2}) by adding a new variable $Y\in \Mb_{m,n}$:
 \al{\label{pb_rw3} &\underset{\substack{X\in \Qb_{m(n+1)}\\ Y\in\Mb_{m,n}}}{ \mathrm{argmin}} \frac{2}{N-n}\ell ( y^N; X)+\frac{2}{N-n} h_{W}(Y)\nn\\ & \hbox{ s.t. }  X_{00}\succ 0,\; \; X\succeq 0\nn\\ & \hspace*{0.7cm}Y=\Dd(X).} 
 The Lagrangian is 
 \al{\Lc &(X,Y,U,Z)=-\log |X_{00}|+\Sp{\Td(\hat R)}{X}\nn\\ &+ \frac{2}{N-n} \bar h_W(Y)-\Sp{U}{X}
 +\Sp{Z}{\Dd(X)-Y}\nn} where $U\in \Qb_{m(n+1)}$ is such that $U\succeq 0$ and $Z\in \Mb_{m,n}$. 
We proceed to perform the unconstrained minimization of $\Lc$ with respect to the primal variables. It is not difficult to see that:
\al{&\underset{ Y}{\min} \,\,\Lc(X,Y,U,Z)\geq \frac{2}{N-n} \ell(y^N; X)+ \Sp{\Td(Z)-U}{X} \nn\\ &\hspace{0.1cm}  +  \underset{Y}{\min}\,\,
\sum_{j=1}^m q_{jj}(Y)\left(\frac{2 \gamma_{jj}}{N-n} -\sum_{k=0}^n  |(Z_k)_{jj}|\right)
\nn\\ &\hspace{0.1cm} + \sum_{j>h} q_{jh}(Y)\left(\frac{2\gamma_{jh}}{N-n} - \sum_{k=0}^n  |(Z_k)_{jh}|+ |(Z_k)_{hj}|\right).\nn
} 
Accordingly, the latter admits minimum if and only if 
\al{\label{cond1}&\sum_{k=0}^n |(Z_k)_{jh}|+|(Z_k)_{hj}|\leq \frac{2\gamma_{jh}}{N-n},\; \; j>h\nn\\
 & \sum_{k=0}^n |(Z_k)_{jj}|\leq \frac{2\gamma_{jj}}{N-n},\; \; j=h;} moreover, in such a situation  it takes value equal to zero. Thus, \al{\underset{ Y}{\min} \,\,\Lc(X,Y,U,Z)=  -\log |X_{00}|+\langle\Td(\hat R)+\Td(Z)-U,X\rangle\nn} if (\ref{cond1}) holds, otherwise it takes $-\infty$.
Regarding the minimization of $\Lc$ with respect to $X$: the minimum with respect to $X_{00}$ is $X_{00}=(\Td(\hat R)+\Td(Z)-U)_{00}^{-1}$; regarding the other blocks $X_{jh}$, it is not difficult to see that $\Lc$ is bounded below if and only if 
\al{\label{cond2} (\Td(\hat R)+\Td(Z)-U)_{jh}=0,\;\; (j,h)\neq (0,0).} We conclude that:
\al{\underset{ X,Y}{\min} \,\,\Lc(X,Y,U,Z)= \log |(\Td(\hat R)+\Td(Z)-U)_{00}|+m\nn} if (\ref{cond1})-(\ref{cond2}) hold, otherwise it takes $-\infty$.
Substituting $U$ with the new variable $W:=(\Td(\hat R)+\Td(Z)-U)_{00}$, we obtain (\ref{dual_rw2}). Regarding the existence of the solution, notice that Problem (\ref{dual_rw2}) consists in maximizing $\log |W|$ subject to the constraints 
$W\succ0$, $\Td(\hat R)+\Td(Z)\succeq  \left[\begin{array}{cc}W &0 \\0 & 0\end{array}\right]$, and $Z\in \mathbf C$ where $\mathbf C:=\{\, Z  \in\Mb_{m,n} \hbox{ s.t. } (\ref{cond1}) \hbox{ holds} \,\}$
which is a closed convex subset of $\{ Z  \in\Mb_{m,n} \hbox{ s.t. } \tr( Z_0)$ $ \hbox{ bounded}\}$.  Then, the existence of the solution is guaranteed by applying a modified version of Lemma A.1 in \cite{LATENTG}. More precisely, it is not difficult to see that the conclusion of such Lemma does not change by replacing the condition $\tr(Z_0)=0$ with $\tr(Z_0)$ bounded. The remaining part of the proof follows the same lines of the one of Proposition 3.3 in \cite{LATENTG}.
\qed

\subsection*{Proof of Proposition \ref{convMM} and Corollary \ref{cor_conv_RWSL}}
We need the following result.
\lemma \label{razaviyayn2013unified_lemma}\citep[Corollary 1]{razaviyayn2013unified}  Consider the problem 
\al{\label{opt_f}u= & \underset{u\in\Cc}{ \mathrm{argmin}} \; \; f(u)} 
which is assumed to be feasible, i.e. there exists $\bar u \in\Cc$ such that $f(\bar u)$ takes finite value, and $\Cc$ is a convex set.
Let $\Mc$ denote the set of stationary points for (\ref{opt_f}).
 Consider the corresponding MM algorithm
\al{u^{(l+1)}= & \underset{u\in\Cc}{ \mathrm{argmin}} \; \; g(u; u^{(l)})\nn} 
where $g(u;v)$ is such that $g(u; u^{(l)})\geq f(u^{(l)})$, $g(u^{(l)}; u^{(l)})=f(u^{(l)})$, $\nabla_u g(u; u)=\nabla_u f(u)$ and  $g(u;v)$ is continuous in $(u,v)$. If the level set $\Sc:=\{ u\in \Cc \hbox{ s.t. } f(u)\leq f(\bar u)\}$ is compact, then the sequence $u^{(l)}$ converges to the set $\Mc$, i.e.
$\lim_{l\rightarrow \infty}\, d(x^{(l)},\Mc)=0$ where $d(x^{(l)},\Mc):=\inf_{s\in\Mc} \|x-s\|$.
\elemma

We proceed to prove Proposition \ref{convMM} by applying Lemma \ref{razaviyayn2013unified_lemma} to Problem (\ref{pb_MM2}) and its MM algorithm (\ref{pb_MM3uno}). It is not difficult to see all the hypotheses on the approximating function $g$ are satisfied. It remains to prove the feasibility and compactness of the level set. Notice that (\ref{pb_MM2}) and (\ref{pb_MM}) are equivalent, thus we can consider Problem (\ref{pb_MM}). In our case, $\Cc=\{\, X\, \hbox{ s.t. } \;  X_{00} \succ0, \; X\succeq 0\,\}$
which is a convex set, and $f(X)=\sum_{j\geq h } \log(q_{jh}(\Dd(X))+\varepsilon_S)+\ell(y^N;X)$. Problem (\ref{pb_MM}) is feasible indeed it is sufficient to pick $\bar X=I$. 
Consider the level set $\Sc:=\{ \,X\in\Cc \; \hbox{ s.t. }\;  
 \;  f(X) \leq f(\bar X)	\,\}$.
We show that $\Sc$ is closed and bounded, since it is a subset of a finite dimensional vector space, then $\Sc $ is compact. 
Consider a convergent sequence $X^{(l)}\in\Cc$, $l\in\Ns$, such that $X_{00}^{(l)}$ converges to a singular matrix, then $-\log|X_{00}^{(l)}|$ tends to infinity. Then, $f(X^{(l)})\geq m(m+1)/2 \log \varepsilon_S-\log|X_{00}^{(l)}|\rightarrow \infty$. Accordingly,  such a sequence cannot belong to $\Sc$, and thus $\Sc$ is a closed set. We proceed to show that unbounded sequences cannot belong to $\Sc$, i.e. $\Sc$ is bounded. Consider a convergent sequence $X^{(l)}\in\Cc$, $l\in\Ns$ such that $\|X^{(l)}\|\rightarrow \infty$. This means that $X^{(l)}$ has at least one eigenvalues tending to infinity as $l\rightarrow \infty$. In $f(X^{(l)})$ the term $\tr(\Td(\hat R)X^{(l)})\rightarrow \infty$ because $\Td(\hat R)\succ0$. Moreover, it dominates the other two terms.  
We conclude that $f(X^{(l)}) \rightarrow \infty$, and thus the sequence cannot belong to $\Sc$. We conclude that $\Sc$ is bounded, and thus compact. Thus, we can apply Corollary 1 and hence the claim is proved. Finally, Corollary \ref{cor_conv_RWSL} can be proved in a similar way.
\qed

\subsection*{Proof of Proposition \ref{prop_norm_const}}
Before to prove the statement we need the following Lemma.
\lemma \label{lemma_integrale}Let $\gamma>0$. Consider the integral: 
\al{\Ic= \int_0^\infty \ldots  \int_0^\infty e^{-\gamma \max \{|x_1|, \ldots ,|x_n|\}}\mathrm d x_1  \ldots \mathrm d x_n.\nn} Then,
 $ \gamma^{-n}\leq \Ic \leq n^n\gamma^{-n}.$\elemma
\proof Notice that:
\al{\frac{1}{n} \sum_{k=1}^n |x_k|\leq \max \{|x_1|, \ldots ,|x_n|\} \leq \sum_{k=1}^n |x_k|,\nn}
accordingly
\al{e^{-\gamma \sum_{k=1}^n |x_k|} \leq e^{-\gamma \max \{|x_1|, \ldots ,|x_n|\} }\leq e^{-\gamma \frac{1}{n} \sum_{k=1}^n |x_k|}.\nn}
Taking the integral with respect to $x_1,\ldots, x_k$ over $[0,\infty)\times\ldots \times [0,\infty)$ in the above inequalities, we obtain a lower and an upper bound for $\Ic$, respectively:
\al{\Ic_L&=\int_0^\infty \ldots \int_0^\infty  e^{-\gamma \sum_{k=1}^n |x_k|} \mathrm d x_1 \ldots \mathrm d x_n \nn\\ &=\prod_{k=1}^n \int_0^\infty e^{-\gamma  |x_k|} \mathrm d x_k= \prod_{k=1}^n \gamma^{-1}=\gamma^{-n}\nn \\
\Ic_U &= \int_0^\infty \ldots \int_0^\infty e^{-\gamma \frac{1}{n} \sum_{k=1}^n |x_k|}\mathrm d x_1 \ldots \mathrm d x_n \nn\\  & = \prod_{k=1}^n \int_0^\infty e^{-\gamma \frac{ |x_k|}{n}} \mathrm d x_k=\prod_{k=1}^n n \gamma^{-1}=n^n \gamma^{-n}\nn } 
which concludes the proof.
\qed

We proceed to prove Proposition \ref{prop_norm_const}.  Since $c_{jh}$ is the normalizing constant of $p_{\gamma_{jh}}((\Sigma)_{jh})$, then we have that 
\al{c_{jh}=\int_{\Qc_{m,n}^+} e^{-\bar h_W(\Sigma)}\mathrm d \Sigma\leq \tilde c_{jh}} where 
\al{\tilde c_{jh}&=\int_{\Qc_{m,n}} e^{-\bar h_W(\Sigma)}\mathrm d \Sigma\nn\\
&=\int_{-\infty}^\infty \ldots \int_{-\infty}^\infty \tilde q_{jh}(S) \mathrm d (S_0)_{jh} \ldots \mathrm d (S_n)_{jh} \mathrm d (S_n)_{hj}\nn} and
\al{\tilde q_{jh}(S)=e^{-\gamma_{jh}\max\{ |(S_0)_{jh}|, \max_{k>0} |(S_k)_{jh}|, \max_{k>0} |(S_k)_{hj}|\}}.\nn} In the case that $j\neq h$, $\tilde c_{jh}$ is an integral taken with respect to $2n+1$ variables.
Since $\tilde q_{jh}(S)=\tilde q_{jh}(-S)$, then  
\al{\tilde c_{jh}=2^{2n+1}\int_0^\infty \ldots \int_0^\infty \tilde q_{jh}(S) \mathrm d (S_0)_{jh} \ldots \mathrm d (S_n)_{jh} \mathrm d (S_n)_{hj}.\nn}
 Hence, by Lemma \ref{lemma_integrale} we have that
\al{\label{form_ofd} \tilde c_{jh}\leq 2^{2n+1} (2n+1)^{2n+1} \gamma^{-(2n+1)}.}
In the case that $j=h$, we have that $(S_k)_{jh}=(S_k)_{hj}$ for $k\geq 0$, accordingly
\al{\tilde c_{jh}=\int_{-\infty}^\infty \ldots \int_{-\infty}^\infty \tilde q_{jh}(S) \mathrm d (S_0)_{jh} \ldots \mathrm d (S_n)_{jh}\nn}
and the integration is taken with respect to $n+1$ variables. Along the same reasoning of before we have that 
\al{\label{form_dia}  \tilde c_{jh}\leq 2^{n+1} (n+1)^{n+1} \gamma^{-(n+1)}.}
In view of (\ref{form_ofd}) and (\ref{form_dia}), we get the claim where $v_{jh}=2^{2n+1}(2n+1)^{2n+1}$ or  $v_{jh}=2^{ n+1}( n+1)^{ n+1}$.
\qed

\subsection*{Proof of Proposition  \ref{prop_opt_gamma}}
Consider the case $j=h$, we have to minimize the objective function $f_{jh}(\gamma_{jh}):=\gamma_{jh} q_{jh}(\hat S^{(l)})-(n+1)\log \gamma_{jh}$
under the constraint that $\gamma_{jh}> 0$. Notice that $f_{jh}(\gamma_{jh})$ is strictly convex. Thus, the minimum point is given by setting equal to zero the
first derivate: $ q_{jh}(\hat S^{(l)})-(n+1)\gamma_{jh}^{-1}=0$ which is satisfied with $\gamma_{jh}=(n+1)/q_{jh}(\hat S^{(l)})$ and the constraint $\gamma_{jh}>0$ is satisfied. The proof for the case $j>h$ is similar. 
\qed

\subsection*{Proof of Theorem \ref{teo_SL_exist}}
The proof is similar to the one of Theorem \ref{teo_S1}. \qed

\subsection*{Proof of Proposition \ref{prop_opt_Q}}
The objective function in (\ref{opt_Q}) is strictly convex in $Q$. Accordingly, the point of minimum is given by setting equal to zero the first variation of the objective function for any $\delta Q\in \Qb_m$:
$\tr (\hat L_0^{(l)}\delta Q-\frac{m+1}{2}Q^{-1}\delta Q )=0$, which implies that 
$ \hat L_0^{(l)}-\frac{m+1}{2} Q^{-1}=0$. Therefore,
under the assumption that $\hat L_0^{(l)}\succ0$ we have (\ref{opt_Q2}).\qed

\end{document}